\documentclass{article}
\usepackage{a4wide}
\usepackage[utf8]{inputenc}
\usepackage{stmaryrd,mathrsfs,bm,amsthm,mathtools,yfonts,amssymb}
\usepackage{color}

\newtheorem{Theorem}{Theorem}[section]
\newtheorem{Proposition}{Proposition}[section]
\newtheorem{Lemma}{Lemma}[section]
\newtheorem{Corollary}{Corollary}[section]
\newtheorem{Definition}{Definition}[section]

\newcommand{\bTheorem}[1]{
\begin{Theorem} \label{T#1} }
\newcommand{\eT}{\end{Theorem}}

\newcommand{\bProposition}[1]{
\begin{Proposition} \label{P#1}}
\newcommand{\eP}{\end{Proposition}}

\newcommand{\bLemma}[1]{
\begin{Lemma} \label{L#1} }
\newcommand{\eL}{\end{Lemma}}

\newcommand{\bCorollary}[1]{
\begin{Corollary} \label{C#1} }
\newcommand{\eC}{\end{Corollary}}

\newcommand{\bFormula}[1]{
\begin{equation} \label{#1}}
\newcommand{\eF}{\end{equation}}

\newcommand{\pat}{\partial_t}
\newcommand{\diver}{\operatorname{div}}

\newcommand{\vr}{\varrho}

\newcommand{\aleq}{\stackrel{<}{\sim}}
\newcommand{\ageq}{\stackrel{>}{\sim}}

\newcommand{\bfphi}{\boldsymbol{\varphi}}

\title{Low stratification of the complete Euler system}

\author{Jan B\v rezina\footnote{Kyushu University, 744 Motooka, Nishi-ku, Fukuoka, 819-0395, Japan; E-mail: brezina@artsci.kyushu-u.ac.jp}, V\' aclav M\' acha\footnote{Institute of Mathematics CAS, \v Zitn\'a 25, 115 67 Praha 1, Czechia; E-mail: macha@math.cas.cz; The work of V\'aclav M\'acha was supported by the Czech Science Foundation project no. GA18-05974S in the framework of RVO: 67985840.}
}

\begin{document}
\maketitle
\abstract{The present paper takes  advantage of the concept of dissipative measure-valued solutions to show the rigorous derivation of the Euler-Boussinesq (EB) system that has been successfully used in various meteorological models. In particular, we show that EB system can be obtained as a singular limit of the complete Euler system. We provide two types of result -- firstly, we treat the case of well-prepared initial  data for any sufficiently regular bounded domain. Secondly,  we use the dispersive estimates for acoustic equation to tackle the case of the ill-prepared initial data on an unbounded exterior domain.}

\section{Introduction}
We consider the flow of an inviscid compressible fluid  governed by the complete  Euler system 
\begin{equation} \label{main.sys}
\begin{split}
\pat \varrho + \diver (\varrho {\bf u}) & = 0\\
\pat(\varrho {\bf u}) + \diver(\varrho {\bf u} \otimes {\bf u}) + \frac 1{{\rm Ma}^2}\nabla p(\varrho, \theta) & = \frac{1}{{\rm Fr}^2} \varrho \nabla F\\
\pat\left(\frac 12 \varrho |{\bf u}|^2 + \frac 1{{\rm Ma}^2} \varrho e (\rho, \theta)\right) + \diver \left[\left(\frac 12 \varrho|{\bf u}|^2 + \frac 1{{\rm Ma}^2}\varrho e (\varrho, \theta)\right){\bf u}\right]&\\ + \frac 1{{\rm Ma}^2}\diver (p(\varrho, \theta){\bf u}) &= \frac 1{{\rm Fr}^2}\varrho \nabla F \cdot {\bf u}
\end{split}
\end{equation}
on a time-space domain $(0,T)\times \Omega$, where $T>0$ and $\Omega \subseteq \mathbb R^3$. The unknowns are scalar quantities $\varrho$ and $\theta$ standing for the density and the temperature, and a vector quantity ${\bf u}$ representing the velocity. Further, we assume that the right-hand side $F$, expressing the potential forces, is independent of time, i.e. $F:\Omega \mapsto \mathbb R$. 

The quantities ${\rm Ma}$ and ${\rm Fr}$ denote Mach  and Froude numbers, respectively. In this paper, we investigate the {\it  low stratification regime}, i.e. we assume that ${\rm Ma}= \varepsilon$ and ${\rm Fr}=\sqrt\varepsilon$ as $\varepsilon \to 0$. This leads to the Euler-Boussinesq system which is of great importance in modelling flows of air in the atmosphere -- we refer to \cite{pedlosky}.

The low Mach number limit  of a system of  Euler equations has been already studied in several papers. The isentropic system has been discussed in \cite{FeKlMa}, where the authors deal with both well and ill-prepared initial data. They consider the case $F=0$.  The complete system has been treated in \cite{BrFe3} and \cite{FeKlKrMa} -- both papers focus on the {\it strong stratification}, i.e. ${\rm Ma} = {\rm Fr} = \varepsilon $ as $\varepsilon \to 0$, for well-prepared initial data. Unlike in the isentropic case, the target system for the complete system case is not uniquely determined and it depends on the choice of the initial data -- the isothermal case has been studied in \cite{FeKlKrMa} and the isentropic case in \cite{BrFe3}.

Concerning the low-Mach number limit of the Navier-Stokes-Fourier system, we refer to a seminal book  by Feireisl and Novotn\'y \cite{FeNo}. Further, the same authors derived the  Euler-Boussinesq equations starting from the Navier-Stokes-Fourier system in \cite{FeNo2}. Besides the low Mach and Froude numbers they also assume that the  Reynolds and P\'eclet numbers are large. Both results are for a monoatomic gas under the effect of radiation.

In this paper we aim to derive the Euler-Boussinesq equations  from the  complete Euler system  under the state equations of a {\it perfect gas}, i.e. we assume
\begin{equation}
p(\varrho,\theta) = \varrho \theta,\qquad e(\varrho, \theta) = c_v \theta,
\end{equation}
where $c_v >0$ is the specific heat at constant volume.
Since the existence of weak solutions to the complete Euler system is currently a difficult task, we take  advantage of a more general concept of  dissipative measure-valued solutions whose existence and properties were discussed recently by B\v rezina and Feireisl (\cite{BrFe}, \cite{BrFe2}, \cite{BrFe4}). In particular, our result holds for finite energy weak solutions as well. 

% is either a regular bounded domain supplemented with the impermeability boundary condition 

Depending on the type of initial data we deal with the following two cases of domains. Either $\Omega\subset \mathbb R^3$ is a  bounded domain or $\Omega\subset\mathbb R^3$ is an exterior domain with compact and smooth boundary -- the regularity of this domain is motivated by an existence of dispersive estimates, see \eqref{dispersive.1} and \eqref{dispersive.2}. For more details we refer to \cite[Section 1 \& 2]{isozaki}. In both situations we assume the impermeability boundary condition 
\begin{equation} \label{bc.imp}
{\bf u} \cdot {\bf n}|_{\partial \Omega} =0.
\end{equation}
In the case of an exterior domain, we  further impose the far field boundary conditions, i.e. 
\begin{equation}\label{far.field}
\varrho\to \overline\varrho, \ \theta\to \overline\theta\mbox{ and }  {\bf u}\to {\bf 0}\ \mbox{as}\ |x|\to \infty,
\end{equation}
where $\overline \varrho$ and $\overline\theta$ are some positive constants specified later. 

As for the potential force $F$ we consider the following situations
\begin{equation}\label{F.specification}
\begin{split}
\Omega \mbox{ bounded: } &F \in W^{k,2}(\Omega) \mbox{ for some }k>5/2, \\
\Omega \mbox{ exterior: } &F(x) = \int_{\mathbb R^3} \frac 1{x-y} m(y)\ {\rm d}y \mbox{ for some }m\mbox{ with}\ {\rm supp}\ m\subseteq \mathbb R^3\setminus\Omega. 
\end{split}
\end{equation}
Note that in the case of an exterior domain the relation in \eqref{F.specification} expresses the fact that the flow of the fluid is driven by the gravitational force of an object lying outside of the fluid.
%\subsection{Entropy production}

Following the Gibbs law, the entropy for the perfect gas is given as
$$
s(\varrho, \theta) = \log\left(\frac{\theta^{c_v}}{\varrho}\right)
$$
and we consider the transport equation for the entropy production rate  in a regularized form
\begin{equation}\label{entropy.production}
\partial_t\left(\varrho \chi\left(\log\left(\frac{\theta^{c_v}}{\varrho}\right)\right)\right) + \diver \left(\varrho \chi \left(\log\left(\frac{\theta^{c_v}}{\varrho}\right)\right){\bf u}\right)\geq 0
\end{equation}
for a certain set of cut-off functions $\chi$ specified later.

Finally, we assume that the total energy of the system is dissipated, i.e.
\begin{equation}\label{energy.conservation}
\frac{{\rm d}}{{\rm d}t} \int_\Omega \left[ \frac 12 \varrho |{\bf u}|^2 + \frac1{{\rm Ma}^2}\varrho e(\varrho, \theta)\right]  \ {\rm d} x \leq \int_\Omega \frac1{{\rm Fr}^2} \varrho \nabla F\cdot {\bf u} \ {\rm d} x.
\end{equation}

On the level of weak solutions we prefer to work  with the system composed of \eqref{main.sys}$_1$, \eqref{main.sys}$_2$ and inequalities \eqref{entropy.production} and \eqref{energy.conservation} rather than \eqref{main.sys}. For one,  the strong solutions to the  systems \eqref{main.sys} and \eqref{main.sys}$_{1,2}$, \eqref{entropy.production}, \eqref{energy.conservation} are equivalent. We refer to \cite{poul} for more details. For two, the latter system allows us to deduce a {\it relative energy inequality} for measure-valued solutions, which is a cornerstone of our proofs.

\subsection{Low stratification and the target system}
In this paper we study exclusively the low stratification regime of the complete Euler system \eqref{main.sys},  i.e. we assume that ${\rm Ma} = \varepsilon$ and ${\rm Fr} = \sqrt{\varepsilon}$ as $\varepsilon \to 0$. We present results for both well and ill-prepared initial data. First, let us formally derive from \eqref{main.sys} the expected limit system  by assuming that
\begin{equation*}
\begin{split}
\varrho &= \overline \varrho + \varepsilon \varrho^{(1)} + \varepsilon^2 \varrho^{(2)} + \ldots,\\
{\bf u} & = {\bf U} + \varepsilon {\bf u}^{(1)} + \varepsilon^2 {\bf u}^{(2)}  + \ldots,\\
\theta & = \overline \theta + \varepsilon \theta^{(1)} + \varepsilon^2 \theta^{(2)} + \ldots.
\end{split}
\end{equation*}
We then multiply \eqref{main.sys}$_2$ by $\varepsilon^2$ and by neglecting all the terms with a positive power of $\varepsilon$ we get
\begin{equation}\label{tlak.limit}
\nabla p(\overline\varrho, \overline\theta)= 0.
\end{equation}
From now on, we  assume that $\overline \theta$  is constant and, consequently, $\overline \varrho$ is also constant due to \eqref{tlak.limit}. These constants agree with the constants presented in \eqref{far.field} for the case of an exterior domain. 

Let us note that the strong stratification limit is quite different as one has to deal with various limit systems (for example isentropic or isothermic) depending on the choice of the relation between $\varrho$ and $\theta$ -- see \cite{BrFe3}. However, in the low stratification case the assumption that the temperature depends on the density  leads immediately to the aforementioned constant states.

It then follows from \eqref{main.sys}$_1$ that
\begin{equation}%\label{isothermal.1}
\diver {\bf U} = 0.
\end{equation}

Further, if we multiply \eqref{main.sys}$_2$ by $\varepsilon$ and omit all the terms with a positive power of $\varepsilon$ we  get the well-known Boussinesq relation
\begin{equation*}%\label{preboussinesq}
\nabla\left(\varrho^{(1)}\overline \theta + \overline \varrho \theta^{(1)}\right) = \overline \varrho \nabla F,
\end{equation*}
which yields the equation
\begin{equation*}
%\label{isothermal.2a}
r +\frac{\overline\varrho}{\overline\theta}\Theta - \frac{\overline\varrho}{\overline\theta}F = h(t),
\end{equation*}
where $r=\varrho^{(1)}$,  $\Theta = \theta^{(1)}$ and  $h(t)$ is some function of time only. Without the loss of generality, it is possible to take $h(t) \equiv 0$ (see \cite[Section 5.3.2]{FeNo}) and hence we use the relation 
\begin{equation}
\label{isothermal.2}
r +\frac{\overline\varrho}{\overline\theta}\Theta = \frac{\overline\varrho}{\overline\theta}F.
\end{equation}

%\comment{\textcolor{blue}{Poh\'adka o $r$: Bu\v d se pou\v zije $r =\varrho^{(1)}$ nebo $r = \varrho^{(1)} - \frac{\overline\varrho}{\overline\theta} F$. Eda d\v el\'a tu druhou verzi, mi p\v rijde logi\v ct\v ej\v s\'i ta prvn\'i. Nicm\'en\v e ob\v e jsou ekvivalentn\'i a p\v ri obou vyjde stejn\'a rovnice \eqref{isothermal.3}, proto\v ze ten cancour nav\'ic, to $\frac{1}{\overline\theta} F\nabla F$, kter\'e by se vyskytlo v \eqref{isothermal.3} pokud bychom uva\v zovali to druh\'e $r$, je vlastn\v e $\nabla\left(\frac 1 {2\overline\theta} F^2\right)$, tedy gradient, kter\'y se d\'a schovat do $\Pi$. }{}{Diky za peknou pohadku. To se ale pak musi nekde pozdeji projevit v REI ne? ale chapu, ze s tim nebude problem tak jako s $\nabla \Pi$.} \textcolor{blue}{p\v resn\v e tak, prost\v e to v\v et\v sinou zmiz\'i, t\v reba d\'iky \eqref{selenoidal.limit}}}

The balance of entropy 
$$
\varrho\partial_t s(\varrho, \theta) + \varrho{\bf u}\cdot \nabla s(\varrho,\theta) = 0
$$
yields 
$$
\overline \varrho Ds(\overline\varrho,\overline\theta)\partial_t(\varrho^{(1)},\theta^{(1)}) + \overline \varrho {\bf U} Ds(\overline \varrho, \overline \theta)\nabla (\varrho^{(1)},\theta^{(1)}) = 0,
$$
which gives
\begin{equation}\label{pre.obe.bou}
\frac{c_v \overline \varrho}{\overline \theta} \partial_t \theta^{(1)} - \partial_t \varrho^{(1)}+ \frac{c_v \overline \varrho}{\overline \theta}  {\bf U} \nabla \theta^{(1)} - {\bf U} \nabla \varrho^{(1)} = 0.
\end{equation}
Moreover, we deduce from \eqref{isothermal.2} that
$$
\varrho^{(1)} = -\frac{\overline \varrho}{\overline \theta} \theta^{(1)} + \frac{\overline \varrho}{\overline \theta} F
$$
and thus \eqref{pre.obe.bou} becomes
\begin{equation}%\label{isothermal.4}
\partial_t \Theta + {\bf U} \cdot \nabla \Theta= \frac1{1+c_v}{\bf U} \cdot \nabla F.
\end{equation}

Finally, we take into account all the terms from \eqref{main.sys}$_2$ with the zero-th power of $\varepsilon$ to get
\begin{equation}%\label{isothermal.3}
 \partial_t {\bf U} + {\bf U}\cdot \nabla{\bf U} + \nabla \Pi  = \frac{r}{\overline \varrho} \nabla F
\end{equation}
for a certain scalar function $\Pi$. 

To summarize we obtained a version of the Euler-Boussinesq system. Namely, we expect that the limit quantities $r$,  ${\bf U},\ \Pi$ and $\Theta$
satisfy
%The system consisting of \eqref{isothermal.1}, \eqref{isothermal.2}, \eqref{isothermal.4} and \eqref{isothermal.3} is a version of the Euler-Boussinesq system. To sum up the heuristic, the quantities ${\bf U},\ \Pi,\ \Theta,$ and $r$ are expected to fulfill:
\begin{equation}\label{Euler.boussinesq}
\begin{split}
\diver {\bf U}& = 0\\
\partial_t {\bf U} + {\bf U}\cdot \nabla{\bf U} + \nabla \Pi &= \frac r{\overline\varrho}\nabla F\\
\partial_t \Theta + {\bf U}\cdot \nabla\Theta & = \frac 1{1+c_v}{\bf U}\cdot \nabla F\\
r + \frac{\overline \varrho}{\overline \theta} \Theta &= \frac{\overline \varrho}{\overline\theta} F 
\end{split}
\end{equation}
with the impermeability boundary condition
\begin{equation}
{\bf U} \cdot {\bf n} = 0.
\end{equation}
Moreover, the  initial conditions  have  the form
\begin{equation}
{\bf U}(0,\cdot) = {\bf U}_0,\ \diver {\bf U}_0 = 0 ,\  \Theta(0,\cdot) = \Theta_0,
\end{equation}
and $r(0,\cdot)$  is given by \eqref{Euler.boussinesq}$_4$ and $\Theta_0$.

Following the nowadays standard theory of well-posedness for hyperbolic systems, see e.g. Kato \cite{Ka}, we can 
assume that the Euler-Boussinesq system \eqref{Euler.boussinesq} with the initial data
\begin{equation} \label{init.data.EB}
\begin{split}
({\bf U}_0,\Theta_0) \in W^{k,2}(\Omega; \mathbb R^4)&, \ \ \ \|( {\bf U_0},\Theta_0 )\|_{W^{k,2}(\Omega, \mathbb R^4)} \leq C,   \\
 \diver{\bf U}_0 = 0&, \ \ \ {\bf U}_0 \cdot {\bf n}|_{\partial \Omega} = 0, \ k> \frac52,
\end{split}
\end{equation}
has a regular solution $( r, {\bf U}, \Theta)$ defined on a maximal time interval $[0, T_{max}(C))$ with the regularity
$$(r , {\bf U},\Theta) \in C([0, T_{max}); W^{k,2}(\Omega; \mathbb R^5)), \ (\partial_t {\bf U}, \nabla \Pi) \in C([0,T_{max}); W^{k-1,2} (\Omega, \mathbb R^6)),$$
where $F$ was specified in \eqref{F.specification}.

\subsection{Conservative variables}

To avoid any problems stemming from the possibility of vacuum we reformulate the Euler system in the conservative variables:
$$
\varrho,\ {\bf m} = \varrho{\bf u},\ p= \varrho \theta.
$$
Here, $\varrho$ is the density, ${\bf m}$ is the momentum and $p$ is the pressure of the fluid.

The system consisting of \eqref{main.sys}$_1$, \eqref{main.sys}$_2$ and inequalities \eqref{entropy.production} and \eqref{energy.conservation} is then equivalent to

\begin{equation}\label{main.sys.conservative}
\begin{split}
\partial_t\varrho + \diver {\bf m} & = 0\\
\partial_t {\bf m} + \diver \left(\frac{{\bf m}\otimes{\bf m}}{\varrho}\right) + \frac 1{\varepsilon^2} \nabla p & = \frac 1\varepsilon \varrho \nabla F\\
\frac{{\rm d}}{{\rm d}t} \int_\Omega \left[ \frac 12 \frac{|{\bf m}|^2}\varrho + \frac{c_v}{\varepsilon^2} p \right]  \ {\rm d} x &\leq \int_\Omega \frac1{\varepsilon} \nabla F \cdot {\bf m} \ {\rm d} x \\
\partial_t\left(\varrho \chi\left(\log\left(\frac{p}{\varrho^\gamma}\right)\right)\right) + \diver \left( \chi \left(\log\left(\frac{p}{\varrho^\gamma}\right)\right){\bf m}\right)&\geq 0, \ \gamma = 1 + \frac1{c_v}
\end{split}
\end{equation}
for a certain set of cut-off functions $\chi$. Here,   \eqref{bc.imp} changes into ${\bf m}\cdot {\bf n}|_{\partial \Omega} = 0$
and \eqref{far.field} into $\varrho \to \overline \varrho$, $p \to \overline \varrho \overline \theta$ and ${\bf m} \to {\bf 0}$ as $|x| \to 0$.

As the functions $\varrho$ and $p$ can in general touch the singular set $\varrho =0$ and/or $p =0$ we extend the "problematic" nonlinearities in 
\eqref{main.sys.conservative} as follows:
\begin{equation}
(\varrho, {\bf m}) \mapsto \frac12 \frac{|{\bf m}|^2}{\varrho} = \left\{ \begin{array}{cl} \frac12 \frac{|{\bf m}|^2}{\varrho} &\mbox{ for } \varrho> 0, \\ 0 & \mbox{ for } {\bf m =0}, \\ \infty & \mbox{ otherwise,} \end{array}\right.
\end{equation}
and
\begin{equation}
(\varrho, p) \mapsto  \varrho \log \left(\frac{p}{\varrho^\gamma} \right) = \left\{ \begin{array}{cl}  \varrho \log \left(\frac{p}{\varrho^\gamma} \right) &\mbox{ for } \varrho\geq  0,\ p>0, \\ -\infty & \mbox{ for } \varrho >0, \ p=0, \\ 0 & \mbox{ otherwise.} \end{array}\right.
\end{equation}

\subsection{Measure-valued solutions to the primitive system}

We establish the set
$$
\mathcal Q \equiv  \left\{(\varrho, {\bf m}, p)|\ \varrho\geq 0, \ {\bf m}\in \mathbb R^3,\ p\geq 0\right\}.
$$

\begin{Definition} We define a {\em renormalized dissipative measure--valued} (rDMV) {\em solution}  to \eqref{main.sys} with the initial conditions $U_0\in L^\infty_{w^*}(\Omega;\mathcal P(\mathcal Q))$ as a parameterized family of probability measures 
$$
U\in L^\infty_{w^*}((0,T)\times \Omega; \mathcal P(\mathcal Q))
$$
that has the following properties:
\begin{itemize}
\item  $\langle U_{t,x};\varrho\rangle - \overline\varrho \in L^{\infty}(0,T,L^1(\Omega))$ for some $\overline\varrho$ positive and
\begin{multline}\label{weak.continuity}
\int_0^T \int_\Omega \left[\langle U_{t,x};\varrho\rangle \partial_t \varphi + \langle U_{t,x};{\bf m}\rangle \cdot \nabla \varphi \right] \ {\rm d }x{\rm d}t =\int_\Omega \langle U_{t,x};\varrho\rangle \varphi(T,\cdot)\ {\rm d}x\\  - \int_\Omega \langle U_{0,x};\varrho \rangle \varphi(0,\cdot)\ {\rm d}x
\end{multline}
for any $\varphi \in C^\infty_c([0,T]\times \overline \Omega)$; 
\item % {}{ $\left\langle U_{t,x}; \frac{|{\bf m}|^2}{\varrho}\right\rangle \in L^\infty(0,T, L^1)$} 
\begin{multline}\label{weak.momentum}
\int_0^T \int_\Omega \left[ \langle U_{t,x}; {\bf m}\rangle \cdot \partial_t \bfphi + \left\langle U_{t,x}; \frac{{\bf m}\otimes{\bf m}}{\varrho}\right\rangle :\nabla \bfphi + \frac 1{\varepsilon^2}\langle U_{t,x};p\rangle \diver \bfphi \right]\ {\rm d}x{\rm d}t\\
 =\int_\Omega \langle U_{T,x};{\bf m}\rangle \cdot \bfphi(T,\cdot)\ {\rm d}x - \int_\Omega \langle U_{0,x};{\bf m}\rangle \cdot \bfphi (0, \cdot) \ {\rm d}x\\
 - \frac 1\varepsilon\int_0^t\int_{\Omega} \langle U_{t,x}; \varrho \rangle \nabla F \cdot \bfphi \ {\rm d}x{\rm d}t  + \int_0^T\int_{\overline \Omega} \nabla \bfphi: {\rm d}\mu_C
\end{multline}
for any $ \bfphi \in C^\infty_c ([0,T]\times \overline \Omega; \mathbb R^3)$, $\bfphi \cdot {\bf n}|_{\partial \Omega} = 0$, where  $\mu_C$ is a vectorial signed measure on $[0,T]\times \overline \Omega$; 
\item  $\langle U_{t,x};\rho (s(\rho,p) - s(\overline\varrho,\overline \varrho\overline\theta))\rangle\in L^\infty(0,T,L^1(\Omega))$ for some $\overline \theta$ positive and
\begin{multline}\label{renorm.entropy}
\int_0^T\int_\Omega \left[ \langle U_{t,x}; \varrho \chi(s(\varrho, p))\rangle \partial_t \varphi + \langle U_{t,x}; \chi (s(\varrho, p)){\bf m}\rangle\cdot \nabla \varphi\right]\ {\rm d}x{\rm d}t\\
\leq\int_\Omega\langle U_{T,x};\varrho\chi(s(\varrho,p))\rangle \varphi\ {\rm d}x -\int_\Omega \langle U_{0,x}; \varrho \chi (s(\varrho, p))\rangle \varphi  \ {\rm d}x
\end{multline}
for any $\varphi \in C_c^\infty ([0,T]\times \overline \Omega)$, $\varphi\geq 0$, and any increasing concave function $\chi$ satisfying $\chi(s)\leq \chi_\infty$ for all $s\in \mathbb R$;

\item The inequality
\begin{equation}\label{weak.energy}
\int_\Omega \left[\left\langle U_{t,x}; \frac12 \frac{|{\bf m}|^2}{\varrho} + \frac 1{\varepsilon^2}  c_v p\right \rangle \right]_{t=0}^{t=\tau}\ {\rm d}x\leq \frac 1{\varepsilon}\int_0^\tau \int_\Omega \langle U_{t,x}; {\bf m}\rangle \cdot \nabla F \ {\rm d}x{\rm d}t
\end{equation}
holds for a.a. $\tau \in [0,T]$.
\item 
\begin{equation} \label{concentration.bound}
\int_0^\tau \int_{\overline \Omega} |{\rm d}\mu_C| \aleq \int_0^\tau \mathcal D(t) \ {\rm d} t
\end{equation}
for a.a. $\tau \in [0,T]$, where 
\begin{equation}\label{concentration}
\mathcal D(\tau) = \int_\Omega \left[\left \langle U_{0,x}; \frac 12 \frac{|{\bf m}|^2}{\varrho} + \frac 1{\varepsilon^2} c_v p \right \rangle   \right]_{t=\tau}^{t=0}\ {\rm d} x+ \frac 1{\varepsilon}\int_0^\tau \int_\Omega \langle U_{t,x};{\bf m}\rangle \cdot \nabla F\ {\rm d}x {\rm d}t.
\end{equation}
Here and hereafter, the symbol $a \aleq b$ means $a \leq cb$ for a certain constant $c>0$ independent of $\varepsilon > 0$, which may  vary from line to line.
\end{itemize} 
\end{Definition}

We tacitly assume that all quantities appearing in the definition are integrable.

We refer to \cite{brezina} and \cite{BrFe}  for  details about the existence of rDMV solutions.

\subsection{Initial data}

For the complete Euler system \eqref{main.sys.conservative} we consider  the initial data in the form
\begin{equation} \varrho_{0,\varepsilon} = \overline \varrho + \varepsilon \varrho_{0,\varepsilon}^{(1)}, \ \ \ {\bf u}_{0,\varepsilon},  \ \ \ \theta_{0,\varepsilon} = \overline \theta + \varepsilon \theta_{0,\varepsilon}^{(1)},\label{IC1}\end{equation}
where $\overline \varrho$, $\overline \theta$ are positive constants and  
\begin{equation}\label{IC2}
 \|\varrho_{0,\varepsilon}^{(1)}\|_{L^\infty(\Omega)}   + \|{\bf u}_{0,\varepsilon}\|_{L^\infty(\Omega)} + \|\theta_{0,\varepsilon}^{(1)}\|_{L^\infty(\Omega)} \leq c.
\end{equation}
Moreover, we assume that 
\begin{equation}\label{init.well.prepared}
\begin{split}
{\bf u}_{0,\varepsilon} &\to {\bf u}_0 \ \mbox{ strongly in }L^2(\Omega),\\
\varrho_{0,\varepsilon}^{(1)} &\to \varrho_{0}^{(1)} \ \mbox{ strongly in }L^2(\Omega),\\
\theta_{0,\varepsilon}^{(1)} &\to  \theta_0^{(1)}\ \mbox{ strongly in }L^2(\Omega).
\end{split}
\end{equation}

\subsubsection{Well-prepared data}
By well-prepared initial data we mean the case when  $\diver {\bf u}_0 =0$ and  $\varrho^{(1)}_{0}$,  $\theta^{(1)}_{0}$ satisfy the Boussinesq  relation \eqref{Euler.boussinesq}$_{(4)}$, i.e.
\begin{equation}\label{boussinesq.ini}
\varrho^{(1)}_0 + \frac{\overline \varrho}{\overline \theta} \theta^{(1)}_0 = \frac{\overline \varrho}{\overline \theta} F.
\end{equation}

We then take
\begin{equation}\label{well.prepared.ini}
r_0 = \varrho^{(1)}_0, \ \ \ {\bf U}_0 = {\bf u}_0, \ \ \  \Theta_0 = \theta^{(1)}_0 \end{equation}
as the initial data for the Euler-Boussinesq system \eqref{Euler.boussinesq}.

The case of well-prepared initial data is treated in  Section \ref{well.prepared.section}.

\subsubsection{Ill-prepared data}

The ill-prepared initial data refer to the general case when neither $\diver {\bf u}_0 =0$ nor  \eqref{boussinesq.ini} need to be satisfied. 
In order to tackle this problem one has to take into account the acoustic analogy. It can be shown that  this discrepancy between the initial data  decays very quickly to $0$ on the exterior domain -- see \eqref{acoustic.equation} and \eqref{dispersive.2}.
In this case we take 
\begin{equation}\label{ill.prepared.ini}
r_0 = \varrho_0^{(1)}, \ \ \ {\bf U}_0 = {\bf H}[{\bf u}_0], \ \ \ \Theta_0 = \frac{\overline \theta}{c_v + 1} \left(c_v\frac1{\overline \theta} \theta_0^{(1)} - \frac 1 {\overline \varrho} \varrho_0^{(1)} + \frac 1{\overline \theta}F\right)
\end{equation}
 as the initial data for \eqref{Euler.boussinesq}. Here, ${\bf H}$ denotes the standard {\it Helmholtz projection} on the space of solenoidal functions. 

The case of ill-prepared initial data is treated in  Section \ref{ill.prepared.section}.

\section{Main results}
We begin with a result for the case of  well-prepared initial data.

\begin{Theorem}
Let $\Omega$ be a regular bounded domain in $\mathbb R^3$, $F$ satisfy \eqref{F.specification} and $\overline \varrho, \overline \theta$ be positive constants. Assume that there is a smooth solution $r$, ${\bf U}$ and $\Theta$ to the Euler-Boussinesq system \eqref{Euler.boussinesq} on $[0,T_{max})$, $T_{max} > 0$ emanating from the initial data
$$
{\bf U_0},\ \ \  \Theta_0
$$
satisfying \eqref{init.data.EB}. % and $r$ given by  $\Theta$ through the relation \eqref{Euler.boussinesq}$_{(4)}$. 
Assume further that $U_{t,x}^\varepsilon$ is a family of  rDMV solutions to \eqref{main.sys} with the initial data 
$$U_{0,x}^\varepsilon = \delta_{\varrho_{0,\varepsilon}, (\varrho_{0,\varepsilon} {\bf u}_{0,\varepsilon}), (\varrho_{0,\varepsilon} \theta_{0,\varepsilon})}, $$
where $ (\varrho_{0,\varepsilon}, {\bf u}_{0,\varepsilon}, \theta_{0,\varepsilon})$ satisfy  
$$
\varrho_{0,\varepsilon} >0,\ \ \ \theta_{0,\varepsilon} > 0,\ \ \    \log\left(\frac{\theta_{0,\varepsilon}^{c_v}}{\varrho_{0,\varepsilon}} \right) \geq s_0 > -\infty
$$ together with \eqref{IC1}--\eqref{well.prepared.ini}.
Moreover, suppose that the constant in \eqref{concentration.bound} is independent of $\{U^\varepsilon\}$.

Then for any $0< T < T_{max}$ there holds  $\mathcal D^\varepsilon \to 0$ in $L^\infty(0,T)$  and
\begin{equation*}
\begin{split}
\langle U_{t,x}^\varepsilon; \varrho\rangle \to \overline\varrho & \ \mbox{strongly in}\ L^\infty(0,T, L^1(\Omega)),\\
\langle U_{t,x}^\varepsilon; p\rangle \to \overline\varrho\overline\theta & \ \mbox{strongly in}\ L^\infty(0,T,L^1(\Omega)),\\
\left\langle U_{t,x}^\varepsilon; \frac{\bf m}{\sqrt \varrho}\right\rangle \to \sqrt{\overline\varrho} {\bf U} & \ \mbox{strongly in}\ L^\infty(0,T,L^2(\Omega)).
\end{split}
\end{equation*}
Moreover,
$$
\left\langle U_{t,x}^\varepsilon; \frac{\varrho - \overline\varrho}{\varepsilon}\right\rangle \to r\ \mbox{and}\ \left\langle U_{t,x}^\varepsilon; \frac{p - \overline\varrho\overline \theta}{\varepsilon}\right\rangle \to \overline\theta r + \overline\varrho \Theta = \overline \varrho F
$$
strongly in $L^\infty(0,T, L^1(\Omega))$. 
\end{Theorem}

Next we formulate the result for  ill-prepared initial data. 

\begin{Theorem}
Let $\Omega \subset \mathbb R^3$ be a regular exterior domain, $F$ satisfy \eqref{F.specification} and $\overline \varrho, \overline \theta$ be positive constants. Assume that there is a smooth solution $r$, ${\bf U}$ and $\Theta$ to the Euler-Boussinesq system \eqref{Euler.boussinesq} on $[0,T_{max})$, $T_{max} >0$ emanating from the initial data
$$
{\bf U_0},\ \ \  \Theta_0
$$
satisfying \eqref{init.data.EB}. % and $r$ given by  $\Theta$ through the relation \eqref{Euler.boussinesq}$_{(4)}$. 
Assume further that $U_{t,x}^\varepsilon$ is a family of  rDMV solutions to \eqref{main.sys} with the initial data 
$$U_{0,x}^\varepsilon = \delta_{\varrho_{0,\varepsilon}, (\varrho_{0,\varepsilon} {\bf u}_{0,\varepsilon}), (\varrho_{0,\varepsilon} \theta_{0,\varepsilon})}, $$
where $ (\varrho_{0,\varepsilon}, {\bf u}_{0,\varepsilon}, \theta_{0,\varepsilon})$ satisfy 
$$
\varrho_{0,\varepsilon} >0,\ \ \   \theta_{0,\varepsilon} > 0, \ \ \   \log\left(\frac{\theta_{0,\varepsilon}^{c_v}}{\varrho_{0,\varepsilon}} \right) \geq s_0 > -\infty
$$
together with \eqref{IC1}--\eqref{init.well.prepared} and \eqref{ill.prepared.ini}.
Moreover, suppose that the constant in \eqref{concentration.bound} is independent of $\{U^\varepsilon\}$.

Then for any $0< T < T_{max}$ there holds $\mathcal D^\varepsilon \to 0$ in $L^\infty(0,T)$  and
\begin{equation*}
\begin{split}
\langle U_{t,x}^\varepsilon; \varrho\rangle \to \overline\varrho & \ \mbox{in strongly}\ L^\infty(0,T, L^1(\Omega)),\\
\langle U_{t,x}^\varepsilon; p\rangle \to \overline\varrho\overline\theta & \ \mbox{in strongly}\ L^\infty(0,T,L^1(\Omega)),\\
\left\langle U_{t,x}^\varepsilon; \frac{\bf m}{\sqrt \varrho}\right\rangle \to \sqrt{\overline\varrho} {\bf U} & \ \mbox{in strongly}\ L^\infty_{loc}((0,T],L^2(\Omega)).
\end{split}
\end{equation*}
Moreover,
$$
\left\langle U_{t,x}^\varepsilon; \frac{\varrho - \overline\varrho}{\varepsilon}\right\rangle \to r\ \mbox{and}\ \left\langle U_{t,x}^\varepsilon; \frac{p - \overline\varrho\overline \theta}{\varepsilon}\right\rangle \to \overline\theta r + \overline\varrho \Theta= \overline \varrho F
$$
strongly in $L_{loc}^\infty((0,T], L^1_{loc}(\overline\Omega))$. 
\end{Theorem}

The idea of the presented proofs is to compare $U_{t,x}^\varepsilon$ with the smooth solution $r,\ {\bf U}, \ \Theta$ using the {\it relative energy inequality} \eqref{rel.inequality}. This is the content of the rest of this paper.

\section{Relative energy functional and inequality}

Fix $\varepsilon>0$ and an increasing concave function $\chi$ satisfying $\chi(s)\leq \chi_\infty$ for all $s\in \mathbb R$  and $\chi(s(\overline\rho,\overline p)) = s(\overline\rho,\overline p)$.  Let $\tilde \varrho >0$, $\tilde{\bf U}$ with  $\tilde {\bf U} \cdot {\bf n}|_{\partial \Omega} = 0$ and $\tilde \theta>0$ be a triple of smooth functions such that  $\tilde\rho - \overline\rho$, $\tilde \theta - \overline\theta$ and $|\tilde{\bf U}|$ decay sufficiently fast to $0$ as $|x|\to \infty$. We then define the functional $\mathcal E_\chi^\varepsilon$   through the following relation:
\begin{multline*}
\mathcal E_\chi^\varepsilon \left(\varrho, {\bf m}, p |\tilde \varrho,\tilde{\bf U}, \tilde\theta\right):= \frac 12 \frac{|{\bf m}|^2}\varrho  - {\bf m}\tilde{\bf U} + \frac 12 \varrho |\tilde{\bf U}|^2\\
+ \frac 1{\varepsilon^2} \left(  c_v p -  \varrho \tilde \theta \chi(s(\varrho, p)) +  \tilde p  - c_v \varrho  \tilde \theta + {\varrho \tilde \theta}s(\tilde \varrho,\tilde p) -  \varrho \tilde \theta \right),
\end{multline*} where we use the notation $\tilde p = \tilde \vr \tilde \theta$.

We deduce from \eqref{weak.continuity}--\eqref{concentration} and the Gibbs law the following:
\begin{multline}\label{rel.inequality}
\left[\int_\Omega \langle U^\varepsilon_{t,x};\mathcal E_\chi^\varepsilon(\varrho,{\bf m}, p|\tilde \varrho,\tilde{\bf U},\tilde \theta)\rangle\ {\rm d}x\right]_{t=0}^{t=\tau} + \mathcal D^\varepsilon(\tau)\\
\leq -\frac 1{\varepsilon^2}\int_0^\tau \int_\Omega \langle U^\varepsilon_{t,x};\varrho \chi(s(\varrho, p))\rangle \partial_t \tilde\theta + \langle U_{t,x}^\varepsilon; \chi(s(\varrho, p)){\bf m}\rangle\cdot \nabla \tilde\theta\ {\rm d}x{\rm d}t\\
+\int_0^\tau \int_\Omega \langle U_{t,x}^\varepsilon;\varrho \tilde{\bf U} - {\bf m}\rangle \cdot \partial_t \tilde{\bf U} + \left\langle U_{t,x}^\varepsilon;\frac{(\varrho\tilde{\bf U} - {\bf m})\otimes{\bf m}}{\varrho}\right \rangle :\nabla \tilde{\bf U} \ {\rm d}x{\rm d}t\\
- \frac 1{\varepsilon^2} \int_0^\tau \int_\Omega \langle U_{t,x}^\varepsilon; p\rangle \diver \tilde{\bf U}\ {\rm d}x{\rm d}t\\
+ \frac 1{\varepsilon^2} \int_0^\tau \int_\Omega \langle U_{t,x}^\varepsilon; \varrho \rangle \partial_t \tilde \theta s(\tilde \vr, \tilde p) + \langle U_{t,x}^\varepsilon;{\bf m}\rangle \cdot \nabla \tilde\theta s(\tilde \vr,\tilde p)\ {\rm d}x{\rm d}t\\
+\frac 1{\varepsilon^2}\int_0^\tau \int_\Omega \langle U_{t,x}^\varepsilon; \tilde\varrho-\varrho\rangle \frac 1{\tilde \varrho} \partial_t \tilde p - \langle U_{t,x}^\varepsilon;{\bf m}\rangle \cdot \frac 1{\tilde\varrho} \nabla \tilde p\ {\rm d}x{\rm d}t\\
+\frac 1\varepsilon \int_0^\tau \int_\Omega  \langle U_{t,x}^\varepsilon;{\bf m} - \varrho \tilde{\bf U}\rangle \cdot \nabla F\ {\rm d}x{\rm d}t + \int_0^\tau \int_\Omega \nabla \tilde{\bf U}: {\rm d}\mu_C^\varepsilon,
\end{multline}
for a.a. $\tau \in (0,T)$.

Furthermore, since it is not difficult to see that 
\begin{equation} \label{rel.en.drop.chi}
- \int_\Omega \langle U_{t,x}^\varepsilon; \varrho \chi (s(\varrho, p))\rangle -  \varrho s(\overline\varrho,\overline p) \ {\rm d}x  \ageq - \int_\Omega \langle U_{t,x}^\varepsilon; \varrho  s(\varrho, p)\rangle-   \varrho  s(\overline\varrho,\overline p) \ {\rm d}x,
\end{equation}  
we shall also use the  {\it relative energy }
\begin{multline*}
\mathcal E^\varepsilon \left(\varrho, {\bf m}, p |\tilde \varrho,\tilde{\bf U}, \tilde\theta\right):= \frac 12 \frac{|{\bf m}|^2}\varrho  - {\bf m}\tilde{\bf U} + \frac 12 \varrho |\tilde{\bf U}|^2\\
+ \frac 1{\varepsilon^2} \left(  c_v p - \varrho \tilde \theta  s(\varrho, p) +  \tilde p  - c_v \varrho  \tilde \theta + {\varrho \tilde \theta}s(\tilde \varrho,\tilde p) -  \varrho \tilde \theta \right).
\end{multline*} independent of any renormalization $\chi$.

\section{Preliminary calculations}
 Similarly to \cite[Section 3.2]{BrFe4} we introduce the {\it essential} and {\it residual} parts of a function $G(\varrho, {\bf m}, p)$ which characterize its behavior near and away from the equilibrium state $(\overline \varrho, \overline \varrho \overline \theta)$, respectively. 

First, we introduce a cutoff function
\begin{multline}\Psi \in C_c^\infty((0,\infty)^2), \ \ 0\leq \Phi \leq 1, \ \ \Phi|_{\mathcal U} =1,\\ \mbox{ where } {\mathcal U} \mbox{ is an open neighborhood  of 
}(\overline \varrho, \overline \varrho \overline \theta) \mbox{ in }(0,\infty)^2.\end{multline}
Second, we decompose any measurable function $G(\varrho, {\bf m}, p)$ as
$$G(\varrho, {\bf m}, p) = [G(\varrho, {\bf m}, p) ]_{ess} + [G(\varrho, {\bf m}, p) ]_{res}$$
with 
$$[G(\varrho, {\bf m}, p) ]_{ess} = \Phi(\varrho, p)G(\varrho, {\bf m}, p)  \mbox{ and } [G(\varrho, {\bf m}, p) ]_{res} = (1- \Phi(\varrho, p))G(\varrho, {\bf m}, p) .$$

Whenever $\bigcup_{(t,x) \in [0,T]\times \Omega} (\tilde \varrho(t,x), \tilde \theta(t,x)) \subseteq K \subseteq \mathcal U$ for some compact set $K$, we can use the following estimate of the  relative energy $\mathcal E^\varepsilon$ that can be found in \cite[Section 3.2.2.]{BrFe4}:
\begin{multline}\label{rel.est}
\mathcal E^\varepsilon\left(\varrho, {\bf m}, p | \tilde \varrho, \tilde{\bf U}, \tilde \theta\right) \ageq \left[ \varrho\left|\frac{{\bf m}}{\varrho} - \tilde{\bf U}\right|^2\right] \\ + \frac 1{\varepsilon^2}\left[|\varrho - \tilde \varrho|^2 + |p - \tilde \varrho \tilde \theta|^2\right]_{ess} + \frac 1{\varepsilon^2} \left[1 + \varrho + \varrho|s(\varrho, p)| +p\right]_{res}.
\end{multline}

For $\tilde \varrho = \overline \varrho$, $\tilde \theta = \overline \theta $, and $\tilde {\bf U } = {\bf 0}$ 
the inequality \eqref{rel.inequality} together with \eqref{rel.en.drop.chi} yield
$$
\int_\Omega \langle U_{\tau,x}^\varepsilon; \mathcal E^\varepsilon (\varrho, {\bf m}, p| \overline \varrho, {\bf 0}, \overline \theta) \rangle \ {\rm d}x \aleq \int_\Omega \langle U_{0,x}^\varepsilon; \mathcal E_\chi^\varepsilon (\varrho, {\bf m}, p| \overline \varrho,{\bf  0}, \overline \theta )\rangle \ {\rm d}x
$$
for all $\tau \in (0,T)$.  As the initial data remain bounded for all $\varepsilon>0$ due to \eqref{IC2}, we deduce that 
\begin{equation}\label{rel.est.2}
\int_\Omega \langle U_{\tau,x}^\varepsilon; [\varrho - \overline\varrho]_{ess}^2 + [p - \overline \varrho\overline \theta]_{ess}^2 + 1_{res} + \varrho_{res} +[\varrho |s(\varrho, p)|]_{res} + p_{res}\rangle \ {\rm d}x \aleq \varepsilon^2 
\end{equation}
and 
$$
\int_\Omega \left\langle U_{\tau,x}^\varepsilon;\frac{|{\bf m}|^2}\varrho \right \rangle \ {\rm d}x \aleq 1
$$
for a.a. $\tau \in (0,T)$.
Therefore we get
\begin{equation}\label{rel.est.3}
{\rm ess\ sup}_{t\in [0,T]} \int_\Omega \langle U_{t,x}^\varepsilon; |\varrho - \overline \varrho|\phi_\varepsilon \rangle \ {\rm d}x \to 0\ \ \ \mbox{ and }\ \ \  {\rm ess\ sup}_{t\in [0,T]} \int_\Omega \langle U_{t,x}^\varepsilon; |p - \overline \varrho\overline \theta|\phi_\varepsilon \rangle \ {\rm d}x \to 0 
\end{equation}
as $\varepsilon \to 0$ for any family  $\{\phi_\varepsilon\}$ bounded in $L^2(\Omega)$ uniformly with respect to $\varepsilon$.

Furthermore, since \eqref{weak.continuity} can be rewritten as
\begin{multline}
 \int_\Omega \langle U_{\tau,x}^\varepsilon; \varrho - \overline \varrho \rangle \varphi(\tau,\cdot) \ {\rm d}x - \int_0^\tau \int_\Omega \langle U_{t,x}^\varepsilon; \varrho - \overline \varrho \rangle \partial_t \varphi  + \langle U_{t,x}^\varepsilon; {\bf m} \rangle \cdot \nabla \varphi \ {\rm d}x{\rm d}t\\ = \int_\Omega \langle U_{0,x}^\varepsilon; \varrho - \overline \varrho \rangle \varphi(0,\cdot ) \ {\rm d}x
\end{multline}
we deduce with the help of \eqref{rel.est.3} that 
\begin{equation}\label{solenoidal.limit}
\int_0^\tau\int_\Omega \langle U_{t,x}^\varepsilon;{\bf m}\rangle \cdot\nabla \varphi \ {\rm d}x{\rm d}t \to 0\ \mbox{as}\ \varepsilon \to 0
\end{equation}
for a.a. $\tau \in [0,T]$ and any $\varphi \in C^1([0,T]\times \overline \Omega)$.

Next, it can be easily shown that once the initial entropy is bounded below it remains bounded below for all times $t>0$. Indeed, since we assume
  $s(\varrho_{0,\varepsilon}, \varrho_{0,\varepsilon} \theta_{0,\varepsilon}) \geq s_0 > -\infty$ we consider \begin{equation} \chi(s) =  \left\{ \begin{array}{cr} s- s_0 & \mbox{ for } s < s_0, \\  0 &\mbox{ for } s\geq s_0. \end{array}\right.\end{equation} 
As this particular $\chi$ is Lipschitz, we deduce from \eqref{weak.continuity} and \eqref{renorm.entropy} that
\begin{multline}
\int_0^\tau\int_\Omega \left[ \langle U_{t,x}; \varrho \Big(\chi(s(\varrho, p)) -\chi(s(\overline\varrho, \overline\varrho\overline\theta))\Big)\rangle \partial_t \varphi + \langle U_{t,x}; \Big(\chi (s(\varrho, p)) - \chi(s(\overline\varrho, \overline\varrho\overline\theta))\Big){\bf m}\rangle\cdot \nabla \varphi\right]\ {\rm d}x{\rm d}t\\
\leq\left[\int_\Omega\langle U_{T,x};\varrho\Big(\chi(s(\varrho,p)) - \chi(s(\overline\varrho,\overline\varrho\overline\theta))\Big)\rangle \varphi\ {\rm d}x\right]_{t=0}^{t=\tau},
\end{multline}
where we can take $\varphi (t,x) \equiv 1$ even  for $\Omega$ exterior. Consequently,
$$ \int_\Omega \langle U_{t,x}^\varepsilon; \varrho \chi (s(\varrho, p))\rangle \ {\rm d}x \geq  \int_\Omega \langle U_{0,x}^\varepsilon; \varrho \chi (s(\varrho, p))\rangle \ {\rm d}x = 0$$
and hence 
$$U_{t,x} \big\{ (\varrho, {\bf m}, p) \big|   s(\varrho, p) \geq s_0 |  \varrho >0 \big\} = 1 \mbox{ for a.a. } (t,x) \in (0,T)\times \Omega.$$
In turn, this means that there exists $a\in \mathbb R$ such that 
\begin{equation} \label{chi1}
 \int_\Omega \langle U_{t,x}^\varepsilon; \varrho \chi (s(\varrho, p))\rangle \ {\rm d}x  = 0 \mbox{ whenever } \chi \leq 0 \mbox{ with } \chi (s) = 0 \mbox{ for all } s \geq a.
\end{equation}

Finally,  we introduce $\chi = \chi_{a,b} \in BC(\mathbb R)$, 
$$\chi_{a,b} (s) = \left\{ \begin{array}{cl} a & \mbox{ for } s < a \\
                                                        s & \mbox{ for } s \in [a,b] \\ 
										   b & \mbox{ for } s \geq b \end{array} \right.$$
with $a, b$ finite fixed in such a way that 
\begin{equation} \label{chi2}
[\chi_{a,b} ( s(\varrho, p))]_{ess} = \Phi(\varrho, p) \chi_{a,b} (s(\varrho, p )) = \Phi (\varrho, p) s(\varrho, p) = [s(\varrho, p)]_{ess}.
\end{equation}

\section{Well-prepared initial data} \label{well.prepared.section}
Let $r$, ${\bf U}$ and $\Theta$ be a solution to the Euler-Boussinesq system \eqref{Euler.boussinesq} with the initial conditions \eqref{well.prepared.ini}. We use \eqref{rel.inequality} with $\tilde \varrho = \overline \varrho + \varepsilon r$,  $\tilde{\bf U} = {\bf U}$, and $\tilde \theta = \overline \theta + \varepsilon \Theta$  to compute

\begin{multline*}
\left[\int_{\Omega}\langle U_{t,x}^\varepsilon; \mathcal E^\varepsilon_\chi(\varrho, {\bf m}, p|\overline \varrho + \varepsilon r, {\bf U}, \overline \theta+ \varepsilon \Theta) \rangle \ {\rm d}x\right]_{t=0}^{t=\tau} + \mathcal D^\varepsilon(\tau) \\
\leq -\frac1{\varepsilon} \int_0^\tau \int_{\Omega} \langle U_{t,x}^\varepsilon; \varrho \chi(  s(\varrho,p))\rangle   \partial_t \Theta  + \langle U_{t,x}^\varepsilon; {\bf m} \chi(s(\varrho, p))\rangle\nabla \Theta \ {\rm d}x{\rm d}t\\
+\int_0^\tau \int_\Omega \langle U_{t,x}^\varepsilon; \varrho {\bf U } - {\bf m}\rangle \cdot \partial_t {\bf U} + \left\langle U_{t,x}^\varepsilon; \frac{(\varrho{\bf U} - {\bf m})\otimes {\bf m}}{\varrho}\right \rangle :\nabla{\bf U} \ {\rm d}x{\rm d}t\\
 + \frac 1\varepsilon \int_0^\tau \int_{\Omega}\left(\langle U_{t,x}^\varepsilon;\varrho\rangle \partial_t\Theta  + \langle U_{t,x}^\varepsilon; {\bf m}\rangle\cdot  \nabla \Theta\right) s(\overline \varrho  + \varepsilon r, (\overline \rho  + \varepsilon r)(\overline \theta + \varepsilon \Theta)) \ {\rm d}x{\rm d}t\\
 + \frac 1{\varepsilon^2} \int_0^\tau \int_\Omega \langle U_{t,x}^\varepsilon; \overline \varrho + \varepsilon r - \varrho \rangle \frac 1{\overline \varrho + \varepsilon r} \partial_t \left((\overline \varrho + \varepsilon r)(\overline \theta + \varepsilon \Theta)\right) - \langle U_{t,x}^\varepsilon;{\bf m}\rangle \cdot \frac 1{\overline\varrho + \varepsilon r} \nabla \left((\overline \varrho + \varepsilon r)(\overline \theta + \varepsilon \Theta)\right)\ {\rm d}x{\rm d}t\\
 + \frac 1\varepsilon \int_0^\tau \int_\Omega  \left \langle U_{t,x}^\varepsilon; {\bf m} -  \varrho {\bf U}\right \rangle \cdot \nabla F   \ {\rm d}x{\rm d}t  + \int_0^\tau \int_\Omega \nabla {\bf U}: {\rm d}\mu_C^\varepsilon
\end{multline*}
for a.a. $\tau \in (0,T)$.

We fix  $a,b \in \mathbb R$ so that \eqref{chi1} and \eqref{chi2} are satisfied. From now on we  consider only $\chi = \chi_{a,b}$ fixed. 
 Thanks to \eqref{rel.en.drop.chi} we can substitute ${\mathcal E}^\varepsilon$ for  ${\mathcal E}_\chi^\varepsilon$ at time $t = \tau >0$ and due to \eqref{IC2} and \eqref{init.well.prepared} we get 
\begin{equation}\int_{\Omega}\langle U_{0,x}^\varepsilon; \mathcal E^\varepsilon_\chi(\varrho, {\bf m}, p|\overline \varrho + \varepsilon r, {\bf U}, \overline \theta+ \varepsilon \Theta) \rangle \ {\rm d}x \to 0\mbox{ as } \varepsilon \to 0,\end{equation}
and consequently, this quantity may be included into $\omega(\varepsilon)$ below.
 Moreover, since   \begin{equation}  \int_0^\tau \int_\Omega \nabla {\bf U}: {\rm d}\mu_C^\varepsilon \aleq  \int_0^\tau \mathcal D^\varepsilon(t)\ {\rm d}t\end{equation} uniformly, we are only left to deal with 

\begin{multline*}
\int_{\Omega}\langle U_{\tau,x}^\varepsilon; \mathcal E^\varepsilon(\varrho, {\bf m}, p|\overline \varrho + \varepsilon r, {\bf U}, \overline \theta+ \varepsilon \Theta) \rangle \ {\rm d}x + \mathcal D^\varepsilon(\tau) \\
\aleq -\frac1{\varepsilon} \int_0^\tau \int_{\Omega} \langle U_{t,x}^\varepsilon; \varrho \chi(  s(\varrho,p))\rangle   \partial_t \Theta  + \langle U_{t,x}^\varepsilon; {\bf m} \chi(s(\varrho, p))\rangle\nabla \Theta \ {\rm d}x{\rm d}t\\
+\int_0^\tau \int_\Omega \langle U_{t,x}^\varepsilon; \varrho {\bf U } - {\bf m}\rangle \cdot \partial_t {\bf U} + \left\langle U_{t,x}^\varepsilon; \frac{(\varrho{\bf U} - {\bf m})\otimes {\bf m}}{\varrho}\right \rangle :\nabla{\bf U} \ {\rm d}x{\rm d}t\\
 + \frac 1\varepsilon \int_0^\tau \int_{\Omega}\left(\langle U_{t,x}^\varepsilon;\varrho\rangle \partial_t\Theta  + \langle U_{t,x}^\varepsilon; {\bf m}\rangle\cdot  \nabla \Theta\right) s(\overline \rho  + \varepsilon r, (\overline \rho  + \varepsilon r)(\overline \theta + \varepsilon \Theta)) \ {\rm d}x{\rm d}t\\
 + \frac 1{\varepsilon^2} \int_0^\tau \int_\Omega \langle U_{t,x}^\varepsilon; \overline \varrho + \varepsilon r - \varrho \rangle \frac 1{\overline \varrho + \varepsilon r} \partial_t \left((\overline \varrho + \varepsilon r)(\overline \theta + \varepsilon \Theta)\right) - \langle U_{t,x}^\varepsilon;{\bf m}\rangle \cdot \frac 1{\overline\varrho + \varepsilon r} \nabla \left((\overline \varrho + \varepsilon r)(\overline \theta + \varepsilon \Theta)\right)\ {\rm d}x{\rm d}t\\
 + \frac 1\varepsilon \int_0^\tau \int_\Omega  \left \langle U_{t,x}^\varepsilon; {\bf m} -  \varrho {\bf U}\right \rangle \cdot \nabla F   \ {\rm d}x{\rm d}t  + \int_0^\tau \mathcal D^\varepsilon(t)\ {\rm d}t + \omega(\varepsilon)\\
=:I + II + III + IV + V + \int_0^\tau \mathcal D^\varepsilon(t)\ {\rm d}t + \omega(\varepsilon)
\end{multline*}
for a.a. $\tau \in (0,T)$.

Here and hereafter,  $\omega(\varepsilon)$ will be used to capture any term tending to $0$, i.e. $\omega(\varepsilon)\to 0$ as $\varepsilon \to 0$.
In the forthcoming lines we provide estimates of  the terms $I$ to $ V$.

Due to \eqref{Euler.boussinesq}$_{(2)} $we have
\begin{multline}
II =  \int_0^\tau \int_\Omega- \langle U_{t,x}^\varepsilon;\varrho{\bf U} - {\bf m}\rangle \cdot  \nabla \Pi +\langle U_{t,x}^\varepsilon; \varrho {\bf U} - {\bf m}\rangle \cdot \frac r{\overline\varrho} \nabla F\ {\rm d}x{\rm d}t\\
 + \int_0^\tau\int_\Omega \left \langle U_{t,x}^\varepsilon; \frac{(\varrho {\bf U} - {\bf m})\otimes({\bf m} - \varrho{\bf U})}{\varrho}\right\rangle :\nabla {\bf U} \ {\rm d}x{\rm d}t
\end{multline}
 where, according to \eqref{rel.est}  the second integral can be estimated  by
\begin{equation}\label{caligraphic.E}
\int_0^\tau \int_\Omega \langle U_{t,x}^\varepsilon, \mathcal E^\varepsilon(\varrho, {\bf m}, p|\overline\varrho + \varepsilon r, {\bf  U}, \overline \theta + \varepsilon\Theta\rangle \ {\rm d}x{\rm d}t.
\end{equation} 

We employ \eqref{Euler.boussinesq}$_{(4)}$, \eqref{rel.est.2} and  \eqref{solenoidal.limit}  in order to obtain
\begin{equation}
II \leq \frac 1{\overline \theta}\int_0^\tau \int_{\Omega}\langle U_{t,x}^\varepsilon; {\bf m} - \varrho{\bf U} \rangle \cdot \Theta \nabla F \ {\rm d}x{\rm d}t+ \int_0^\tau \int_\Omega \langle U_{t,x}^\varepsilon, \mathcal E^\varepsilon(\varrho, {\bf m}, p|\overline\varrho + \varepsilon r, {\bf  U}, \overline \theta + \varepsilon\Theta\rangle \ {\rm d}x{\rm d}t + \omega(\varepsilon).
\end{equation} 

Next,
\begin{multline} \label{respes}
I + III = -\frac 1\varepsilon \int_0^\tau\int_{\Omega} \langle U_{t,x}^\varepsilon; \varrho (\chi(s(\varrho, p)) - s(\overline \varrho + \varepsilon r,(\overline \varrho + \varepsilon r)( \overline \theta + \varepsilon \theta)))\rangle (\partial_t \Theta + {\bf U} \cdot \nabla \Theta)\ {\rm d}x{\rm d}t\\
- \frac 1\varepsilon \int_0^\tau \int_\Omega \langle U_{t,x}^\varepsilon;({\bf m} - \varrho{\bf U})(\chi(s(\varrho, p)) - s(\overline \varrho + \varepsilon r,(\overline \varrho + \varepsilon r)(\overline \theta + \varepsilon \Theta)))\rangle\cdot  \nabla \Theta\ {\rm d}x{\rm d}t.
\end{multline}
The second integral can be estimated by \eqref{caligraphic.E} and $\omega(\varepsilon)$. Indeed,  as $\chi$ is bounded we can estimate the residual parts  by the Young inequality and \eqref{rel.est} as
\begin{multline}
-\frac 1\varepsilon \int_0^\tau \int_\Omega \langle U_{t,x}^\varepsilon;({\bf m} - \varrho{\bf U})_{res}(\chi(s(\varrho, p)) - s(\overline \varrho + \varepsilon r,(\overline \varrho + \varepsilon r)(\overline \theta + \varepsilon \Theta)))_{res}\rangle\cdot  \nabla \Theta\ {\rm d}x{\rm d}t\\
\aleq \int_0^\tau \int_{\Omega} \left\langle U_{t,x}^\varepsilon; \varrho \left|\frac{\bf m}{\varrho} - {\bf U}\right|^2\right\rangle \ {\rm d}x{\rm d}t + \frac {1}{\varepsilon^2} \int_0^\tau \int_\Omega \left\langle U_{t,x}^\varepsilon; \varrho_{res}\right\rangle \ {\rm d}x{\rm d}t\\
\aleq \int_0^\tau \int_\Omega \langle U_{t,x}^\varepsilon; \mathcal E^\varepsilon(\varrho, {\bf m}, p| \overline\varrho + \varepsilon r, {\bf U}, \overline\theta + \varepsilon\Theta)\rangle\ {\rm d}x{\rm d}t.
\end{multline}
On the essential part $\chi(s(\varrho, p)) = s(\varrho, p)$,  $s$ is a Lipschitz function and once again using the Young inequality and \eqref{rel.est} we obtain
\begin{multline}
-\frac 1\varepsilon \int_0^\tau \int_\Omega \langle U_{t,x}^\varepsilon;({\bf m} - \varrho{\bf U})_{ess}(s(\varrho, p)_{ess} - s(\overline \varrho + \varepsilon r,(\overline \varrho + \varepsilon r)(\overline \theta + \varepsilon \Theta))_{ess})\rangle\cdot  \nabla \Theta\ {\rm d}x{\rm d}t\\
\aleq \int_0^\tau\int_\Omega \langle U_{t,x}^\varepsilon; \frac{|{\bf m} - \varrho{\bf U}|^2}{\varrho}_{ess}\rangle \ {\rm d}x{\rm d}t + \frac 1{\varepsilon^2}\int_0^\tau \int_\Omega \langle U_{t,x}^\varepsilon; (\varrho - \overline\varrho - \varepsilon r)_{ess}^2\rangle\ {\rm d}x{\rm d}t\\
 + \frac 1{\varepsilon^2}\int_0^\tau \int_\Omega \langle U_{t,x}^\varepsilon; (p - (\overline\varrho + \varepsilon r)(\overline\theta + \varepsilon \Theta))_{ess}^2\rangle\ {\rm d}x{\rm d}t\\
\aleq \int_0^\tau \int_{\Omega} \langle U_{t,x}^\varepsilon; \mathcal E^\varepsilon(\varrho, {\bf m}, p|\overline\varrho + \varepsilon r,  {\bf U}, \overline\theta + \varepsilon \Theta\rangle \ {\rm d}x{\rm d}t.
\end{multline}

To deal with  $IV$ we use \eqref{Euler.boussinesq}$_{(4)}$ in order to claim that
\begin{multline}
\frac1{\varepsilon^2}\int_0^\tau \int_\Omega \langle U_{t,x}^\varepsilon;{\bf m}  \rangle \cdot\frac 1{\overline \varrho + \varepsilon r} \nabla\left((\overline \varrho + \varepsilon r)(\overline \theta + \varepsilon \Theta)\right)\ {\rm d}x{\rm d}t \\
= \frac1{\varepsilon^2}\int_0^\tau \int_\Omega \langle U_{t,x}^\varepsilon;{\bf m}\rangle\cdot \frac1{\overline\varrho + \varepsilon r} \varepsilon^2\nabla(r\Theta)\ {\rm d}x{\rm d}t\\ 
+ \frac1{\varepsilon^2}\int_0^\tau \int_\Omega  \langle U_{t,x}^\varepsilon;{\bf m}\rangle \cdot \frac 1{\overline \varrho + \varepsilon r} \varepsilon\left(\overline\theta \nabla r + \overline\varrho \nabla \Theta\right)\ {\rm d}x{\rm d}t\\
=\int_0^\tau \int_\Omega \langle U_{t,x}^\varepsilon;{\bf m}\rangle \cdot \nabla H\ {\rm d}x{\rm d}t + \frac 1\varepsilon \int_0^\tau \int_\Omega \langle U_{t,x}^\varepsilon;{\bf m}\rangle  \cdot \nabla F\ {\rm d}x{\rm d}t\\
 - \int_0^\tau \int_\Omega \langle U_{t,x}^\varepsilon; {\bf m}\rangle \cdot \frac{r}{\overline\varrho + \varepsilon r} \nabla F\ {\rm d}x{\rm d}t+ \omega(\varepsilon),
\end{multline}
where $H$ is some function of $r$ and $\Theta$. Thus, the first integral on the right-hand side tends to zero due to \eqref{solenoidal.limit} and, consequently, it may be included into $\omega(\varepsilon)$. 

Summarizing our effort we arrive at 
\begin{multline}
\int_{\Omega}\langle U_{\tau,x}^\varepsilon; \mathcal E^\varepsilon(\varrho, {\bf m}, p|\overline \varrho + \varepsilon r, {\bf U}, \overline \theta+ \varepsilon \Theta) \rangle \ {\rm d}x + \mathcal D^\varepsilon(\tau) \\
\aleq \int_0^\tau \int_\Omega \langle U_{t,x}^\varepsilon; \mathcal E^\varepsilon(\varrho, {\bf m}, p|\overline \varrho + \varepsilon r, {\bf U}, \overline \theta + \varepsilon \Theta\rangle \ {\rm d}x{\rm d}t + \int_0^\tau \mathcal D^\varepsilon(t)\ {\rm d}t  + \omega(\varepsilon) \\
% -\frac1{\varepsilon} \int_0^\tau \int_{\Omega} \langle U_{t,x}^\varepsilon; \varrho s(\varrho,\theta)\rangle   \partial_t \Theta  + \langle U_{t,x}^\varepsilon; {\bf m} s(\varrho, \theta)\rangle\nabla \Theta \ {\rm d}x{\rm d}t\\
%\int_0^\tau \int_\Omega \langle U_{t,x}^\varepsilon; \varrho {\bf U } - {\bf m}\rangle \cdot \partial_t {\bf U} + \left\langle U_{t,x}^\varepsilon; \frac{(\varrho{\bf U} - {\bf m})\otimes {\bf m}}{\varrho}\right \rangle :\nabla{\bf U} \ {\rm d}x{\rm d}t\\
% + \frac 1\varepsilon \int_0^\tau \int_{\Omega}\left(\langle U_{t,x}^\varepsilon;\varrho\rangle \partial_t\Theta  + \langle U_{t,x}^\varepsilon; {\bf m}\rangle \nabla \Theta\right) s(\overline \rho  + \varepsilon r, \overline \theta + \varepsilon \Theta) \ {\rm d}x{\rm d}t\\
-\frac 1\varepsilon \int_0^\tau \int_\Omega \langle U_{t,x}^\varepsilon; \varrho(\chi(s(\varrho, p)) - s(\overline  \varrho + \varepsilon r, (\overline  \varrho + \varepsilon r)(\overline \theta + \varepsilon \theta)))\rangle (\partial_t \Theta + {\bf U} \cdot \nabla\Theta) \ {\rm d}x{\rm d}t\\
 + \frac 1{\varepsilon^2} \int_0^\tau \int_\Omega \langle U_{t,x}^\varepsilon; \overline \varrho + \varepsilon r - \varrho \rangle \frac 1{\overline \varrho + \varepsilon r} \partial_t \left((\overline \varrho + \varepsilon r)(\overline \theta + \varepsilon \Theta)\right)\ {\rm d}x {\rm d}t\\% - \langle U_{t,x}^\varepsilon;{\bf m}\rangle \frac 1{\overline\varrho + \varepsilon r} \nabla \left((\overline \varrho + \varepsilon r)(\overline \theta + \varepsilon \Theta)\right)\ {\rm d}x{\rm d}t\\
 + \frac 1\varepsilon \int_0^\tau \int_\Omega   \left \langle U_{t,x}^\varepsilon; {\bf m} - \varrho {\bf U}\right \rangle \cdot \nabla F \left(1 + \frac{\varepsilon \Theta}{\overline \theta}\right) \ {\rm d}x{\rm d}t - \frac 1\varepsilon\int_0^\tau \int_\Omega \langle U_{t,x}^\varepsilon;{\bf m}\rangle \cdot \nabla F \ {\rm d}x{\rm d}t
\\ 
+ \int_0^\tau \int_\Omega  \langle U_{t,x}^\varepsilon;{\bf m}\rangle \cdot \frac r{\overline\varrho} \nabla F \ {\rm d}x{\rm d}t
\end{multline} 
for a.a. $\tau \in (0,T)$.

According to \eqref{Euler.boussinesq}$_{(4)}$ we have 
\begin{equation}
\partial_t(\overline \theta r+\overline\varrho \Theta) = 0
\end{equation}
and thus
\begin{multline}
\frac 1{\varepsilon^2} \int_0^\tau \int_\Omega \langle U_{t,x}^\varepsilon; \overline \varrho + \varepsilon r - \varrho \rangle \frac 1{\overline \varrho + \varepsilon r} \partial_t \left((\overline \varrho + \varepsilon r)(\overline \theta + \varepsilon \Theta)\right)\\
= \frac 1{\varepsilon^2} \int_0^\tau \int_\Omega \langle U_{t,x}^\varepsilon; \overline \varrho + \varepsilon r-\varrho\rangle \frac1{\overline \varrho + \varepsilon r} \varepsilon^2\partial_t(r \Theta)\ {\rm d}x{\rm d}t = \omega(\varepsilon).
\end{multline}

Further, we use \eqref{Euler.boussinesq}$_{(3)}$ to obtain 

\begin{multline}\label{well.posed.almost}
\int_{\Omega}\langle U_{\tau,x}^\varepsilon; \mathcal E^\varepsilon(\varrho, {\bf m},p|\overline \varrho + \varepsilon r, {\bf U}, \overline \theta+ \varepsilon \Theta) \rangle \ {\rm d}x + \mathcal D^\varepsilon(\tau) \\
\aleq \int_0^\tau \int_\Omega \langle U_{t,x}^\varepsilon; \mathcal E^\varepsilon(\varrho, {\bf m}, p |\overline \varrho + \varepsilon r, {\bf U}, \overline \theta + \varepsilon \Theta\rangle \ {\rm d}x{\rm d}t + \int_0^\tau \mathcal D^\varepsilon(t)\ {\rm d}t  + \omega(\varepsilon) \\
-\frac 1\varepsilon \int_0^\tau \int_\Omega \langle U_{t,x}^\varepsilon; \varrho(\chi(s(\varrho,p)) - s(\overline  \varrho + \varepsilon r, (\overline  \varrho + \varepsilon r)(\overline \theta + \varepsilon \theta)))\rangle \frac 1{c_v + 1} {\bf U} \cdot \nabla F \ {\rm d}x{\rm d}t\\
+ \int_0^\tau \int_\Omega \langle U_{t,x}^\varepsilon;{\bf m} - \overline \varrho {\bf U}\rangle \cdot \frac{\Theta}{\overline \theta}\nabla F\ {\rm d}x{\rm d}t - \frac 1 \varepsilon \int_0^\tau \int_\Omega\langle U_{t,x}^\varepsilon; \varrho {\bf U}\rangle \cdot  \nabla F\ {\rm d}x{\rm d}t\\
+ \int_0^\tau \int_\Omega  \langle U_{t,x}^\varepsilon;{\bf m}\rangle \cdot \frac r{\overline\varrho} \nabla F \  {\rm d}x{\rm d}t
\end{multline}
for a.a. $\tau \in (0,T)$.

As
\begin{multline*}
\frac 1\varepsilon \int_0^\tau \int_\Omega \langle U_{t,x}^\varepsilon; \varrho{\bf U}\rangle\cdot \nabla F \ {\rm d}x{\rm d}t= \frac1\varepsilon \int_0^\tau \int_\Omega\langle U_{t,x}^\varepsilon; (\varrho - \overline \varrho) {\bf U}\rangle \cdot \nabla F \ {\rm d}x{\rm d}t\\
=\frac1\varepsilon \int_0^\tau \int_\Omega \langle U_{t,x}^\varepsilon; (\varrho - \overline \varrho - \varepsilon r) \rangle {\bf U} \cdot \nabla F \ {\rm d}x{\rm d}t + \int_0^\tau \int_\Omega \langle U_{t,x}^\varepsilon; \overline \varrho{\bf U}\rangle \cdot \frac r{\overline \varrho} \nabla F \ {\rm d}x{\rm d}t,
\end{multline*}
we use \eqref{Euler.boussinesq}$_{(4)}$. We recall that $F\nabla F = \frac 12 \nabla |F|^2$ and thus the term $ \displaystyle \int_0^\tau\int_\Omega \langle U_{t,x}^\varepsilon,\overline\varrho{\bf U}\rangle \cdot \frac1{\overline\theta} F\nabla F \ {\rm d}x{\rm d}t = 0$ due to \eqref{Euler.boussinesq}$_{(1)}$. Consequently, the last three integrals in \eqref{well.posed.almost} can be rewritten as
\begin{equation} \label{almost.end}
- \frac 1\varepsilon\int_0^\tau \int_\Omega \langle U_{t,x}^\varepsilon; (\varrho - \overline \varrho - \varepsilon r)\rangle {\bf U} \cdot  \nabla F \ {\rm d}x{\rm d}t.
\end{equation}

Finally, we handle the "entropy" integral on the right-hand side of \eqref{well.posed.almost}. First we realize that  the residual part tends to zero as $\varepsilon \to 0$ thanks to the boundedness of $\chi$. Second, we can omit $\chi$ on the essential part and we use Taylor expansion (recall that $s(\varrho, p) = c_v \log p - (c_v + 1)\log \varrho)$ in order to get
\begin{multline}\label{taylor}
\frac{s(\varrho,p) - s(\overline\varrho + \varepsilon r, (\overline\varrho + \varepsilon r)(\overline\theta + \varepsilon\Theta))}{\varepsilon}\\ = Ds(\overline\varrho + \varepsilon r, (\overline \varrho + \varepsilon r)(\overline \theta + \varepsilon \Theta))\left(\frac{\varrho - \overline\varrho - \varepsilon r}{\varepsilon}, \frac{p - (\overline \varrho + \varepsilon r)(\overline\theta + \varepsilon\Theta)}{\varepsilon}\right) + \omega(\varepsilon)\\
=-\frac{c_v + 1}{\varepsilon} \frac{\varrho - \overline\varrho - \varepsilon r}{\overline\varrho + \varepsilon r} + \frac{c_v}\varepsilon \frac{p - (\overline\varrho + \varepsilon r)(\overline\theta + \varepsilon\Theta)}{(\overline \varrho + \varepsilon r)(\overline \theta + \varepsilon \Theta)} + \omega(\varepsilon),
\end{multline}
where $\omega(\varepsilon)$ stands for a function whose appropriate integral tends to $0$ as $\varepsilon \to 0$. We note that the relation \eqref{taylor} is valid only on the essential part. 
Consequently,
\begin{multline}
-\frac 1\varepsilon \int_0^\tau \int_\Omega \langle U_{t,x}^\varepsilon; \varrho (\chi(s(\varrho,p)) - s(\overline  \varrho + \varepsilon r, (\overline  \varrho + \varepsilon r)(\overline \theta + \varepsilon \theta)))\rangle \frac 1{c_v + 1} {\bf U} \cdot\nabla F \ {\rm d}x{\rm d}t\\
\aleq \frac 1\varepsilon\int_0^\tau\int_\Omega \langle U_{t,x}^\varepsilon; [\varrho - \overline\varrho - \varepsilon r]_{ess}\rangle  {\bf U}\cdot \nabla F\ {\rm d}x{\rm d}t\\
-\frac{c_v}{\overline \theta(c_v + 1)}\int_0^\tau \int_\Omega \left\langle U_{t,x}^\varepsilon; \left[ \frac{p - (\overline\varrho +\varepsilon r)(\overline\vartheta + \varepsilon\Theta)}{\varepsilon}\right]_{ess}\right\rangle {\bf U}\cdot \nabla F\ {\rm d}x{\rm d}t\\
+\varepsilon \int_0^\tau \int_\Omega \langle U_{t,x}^\varepsilon; \mathcal E^\varepsilon(\varrho, {\bf m}, p|\overline \varrho + \varepsilon r, {\bf U}, \overline \theta + \varepsilon \Theta)\rangle \ {\rm d}x{\rm d}t + \omega(\varepsilon).
\end{multline} 
Note, that the first term on the right-hand side cancels with the essential part of \eqref{almost.end} and the remaining residual part of \eqref{almost.end} tends to zero as $\varepsilon \to 0$. Thus it remains only to handle 
\begin{equation}\label{pressure.term}
-\frac{c_v}{\overline \theta(c_v + 1)}\int_0^\tau \int_\Omega \left\langle U_{t,x}^\varepsilon; \left[ \frac{p - (\overline\varrho +\varepsilon r)(\overline\vartheta + \varepsilon\Theta)}{\varepsilon} \right]_{ess}\right\rangle {\bf U}\cdot \nabla F\ {\rm d}x{\rm d}t.
\end{equation}

In order to do that, we multiply the momentum equation \eqref{weak.momentum} by $\varepsilon$ to get that
\begin{equation}
\int_0^\tau \int_\Omega \frac1\varepsilon \langle U_{t,x}^\varepsilon; p\rangle \diver \bfphi + \langle U_{t,x}^\varepsilon; \varrho \rangle \nabla F \cdot \bfphi \ {\rm d}x{\rm d}t = \omega(\varepsilon)
\end{equation}
for any $\bfphi \in C^\infty_c ([0,T)\times \overline \Omega; \mathbb R^3)$, $\bfphi \cdot {\bf n}|_{\partial \Omega} = 0$. This 
 can be rewritten as 
\begin{equation}\label{pressure.term2}
\int_0^\tau \int_\Omega \frac1\varepsilon \langle U_{t,x}^\varepsilon; p - \overline\varrho \overline \theta\rangle \diver \bfphi + \langle U_{t,x}^\varepsilon; \varrho-\overline \varrho \rangle \nabla F \cdot \bfphi  + \overline \varrho \nabla F\cdot \bfphi\ {\rm d}x{\rm d}t = \omega(\varepsilon),
\end{equation}
We use  integration by parts together with \eqref{Euler.boussinesq}$_{(4)}$  and, consequently, \eqref{pressure.term2} can be rewritten as
\begin{equation} \label{deal.pressure.term}
\int_0^\tau \int_\Omega \left\langle U_{t,x}^\varepsilon; \frac{p-\varepsilon \overline \varrho \Theta - \varepsilon \overline \theta  r - \overline \varrho \overline \theta}{\varepsilon}\right\rangle \diver \bfphi \ {\rm d}x{\rm d}t = \int_0^\tau \int_\Omega \langle U_{t,x}^\varepsilon; \varrho - \overline \varrho\rangle \nabla F \cdot \bfphi + \omega(\varepsilon),
\end{equation}
where the right-hand side tends to zero due to \eqref{rel.est.2}. Taking $\varphi = F{\bf U}$ in \eqref{deal.pressure.term}  and realizing that the residual parts tend to zero as $\varepsilon \to 0$ this yields that \eqref{pressure.term} can be also included into $\omega(\varepsilon)$.

We put all these estimates into \eqref{well.posed.almost} in order to get
\begin{multline}
\int_{\Omega}\langle U_{\tau,x}^\varepsilon; \mathcal E^\varepsilon(\varrho, {\bf m},p|\overline \varrho + \varepsilon r, {\bf U}, \overline \theta+ \varepsilon \Theta) \rangle \ {\rm d}x + \mathcal D^\varepsilon(\tau) \\
\aleq \int_0^\tau \int_\Omega \langle U_{t,x}^\varepsilon; \mathcal E^\varepsilon(\varrho, {\bf m}, p|\overline \varrho + \varepsilon r, {\bf U}, \overline \theta + \varepsilon \Theta)\rangle \ {\rm d}x{\rm d}t + \int_0^\tau \mathcal D^\varepsilon(t)\ {\rm d}t  + \omega(\varepsilon)
\end{multline}
for a.a. $\tau \in (0,T)$  and the Gr\"onwall inequality concludes the proof. 

\section{Ill-prepared initial data}
\label{ill.prepared.section}
We recall that throughout this section  $\Omega$ is an exterior domain with a nonempty complement.  

Let $r$, ${\bf U}$ and $\Theta$ be a solution to the Euler-Boussinesq system \eqref{Euler.boussinesq} with the initial conditions \eqref{ill.prepared.ini}. 
In order to deal with the ill-prepared initial data, we proceed similarly as in \cite{FeNo2}, i.e. we choose  test functions $\{\tilde \varrho, \tilde{\bf U}, \tilde \theta\}$ in the relative energy inequality \eqref{rel.inequality} in the following way:
\begin{equation}\label{test.function.ill.posed}
\tilde \varrho = \overline\varrho + \varepsilon R_\varepsilon, \quad \tilde{\bf U}= {\bf U} + \nabla \Phi_\varepsilon,\quad \tilde \theta = \overline\theta + \varepsilon T_\varepsilon,
\end{equation}
where the functions $R_\varepsilon$, $T_\varepsilon$ and $\Phi_\varepsilon$ satisfy the acoustic equation
\begin{equation}\label{acoustic.equation}
\begin{split}
\varepsilon \partial_t \left(\frac{\overline\theta}{\overline\rho} R_\varepsilon + T_\varepsilon - F\right) + \frac{(c_v+1)\overline\theta}{c_v}\Delta \Phi_\varepsilon& = 0,\ \nabla \Phi\cdot {\bf n} |_{\partial\Omega}= 0,\\
\varepsilon \partial_t \nabla \Phi_\varepsilon  + \nabla \left(\frac{\overline \theta}{\overline\varrho} R_\varepsilon + T_\varepsilon - F\right) & = 0.
\end{split}
\end{equation}
The initial data are given by 
$$
R_\varepsilon(0,\cdot) = R_0,\ T_\varepsilon(0,\cdot) = T_0,\ \Phi_{\varepsilon}(0,\cdot) = \Phi_0.
$$

Further, $T_\varepsilon$ and $R_\varepsilon$ are supposed to fulfill
\begin{equation}\label{temperature.for.T}
\partial_t\left(\frac{c_v\overline \varrho}{\overline\theta} T_\varepsilon - R_\varepsilon\right) + \tilde{\bf U}_\varepsilon\cdot \nabla \left(\frac{c_v\overline\varrho}{\overline\theta} T_\varepsilon - R_\varepsilon\right) = 0.
\end{equation}

In accordance with \cite{FeNo2} (see also \cite{FeNoSu}) we introduce the following regularization:
\begin{equation}
[v]_\eta = G_\eta(\sqrt{-\Delta_N})[\psi_{1/\eta} v],
\end{equation}
where 
$$
\psi_\eta(x) = \psi(x/\eta); \psi \in C_c^\infty(\mathbb R),\ 0\leq \psi \leq 1,\ \psi = \left\{\begin{matrix} 1\ \mbox{for } |x|\leq 1\\ 0\ \mbox{for }|x|\geq 2\end{matrix}\right.
$$
and
\begin{equation}
\begin{split}
G_\eta \in C_c^\infty (\mathbb R),\ 0\leq G_\eta\leq 1,\ G_\eta(-z) = G_\eta(z),\\
G_\eta(z) = 1\ \mbox{for } z\in (\eta,1/\eta),\\
G_\eta(z) = 0\ \mbox{for } z\in [0,\eta/2)\cup (2/\eta,\infty).
\end{split}
\end{equation}
We take  $R_0= R_{0,\eta}$, $T_0= T_{0,\eta}$ and $\Phi_0 =\Phi_{0,\eta}$ as the initial data for \eqref{acoustic.equation} such that:
\begin{equation*}
\begin{split}
\|R_{0,\eta}\|_{\infty}\leq c(\eta),\ R_{0,\eta}& \to \varrho_0^{(1)}\ \mbox{strongly in}\ L^2(\Omega),\\
\|T_{0,\eta}\|_\infty \leq c(\eta), \ T_{0,\eta}& \to \theta_0^{(1)}\ \mbox{strongly in}\ L^2(\Omega),\\
\nabla\Phi_{0,\eta}& \to {\bf H}^\bot[{\bf u}_0]\ \mbox{strongly in}\ L^2(\Omega)
\end{split}
\end{equation*}
as $\eta \to 0$. Here,  ${\bf H}^\perp$ stands for the orthogonal complement of the Helmholtz projection ${\bf H}$ in $L^2(\Omega)$. Consequently, the test functions introduced in \eqref{test.function.ill.posed} depend on $\varepsilon$ and $\eta$ so we adopt the notation $\tilde \varrho_{\varepsilon,\eta}$, $\tilde \theta_{\varepsilon,\eta}$ and $\tilde{\bf U}_{\varepsilon,\eta}$. 

Immediately, due to the choice of the initial data, we have
\begin{equation}\label{rel.ini.ill}
\int_\Omega \langle U^\varepsilon_{0,x};\mathcal E_\chi^\varepsilon(\varrho,{\bf m}, p|\tilde \varrho_{\varepsilon,\eta},\tilde{\bf U}_{\varepsilon,\eta},\tilde \theta_{\varepsilon,\eta})\rangle\ {\rm d}x = \omega(\varepsilon, \eta),
\end{equation}
where $\omega:=\omega(\varepsilon,\eta)$ denotes any quantity that tends to zero, namely
$$
\lim_{\eta\to 0}(\lim_{\varepsilon \to 0} \omega(\varepsilon, \eta)) = 0.
$$

The following set of estimates as well as further details about the equations \eqref{acoustic.equation} and \eqref{temperature.for.T} can be found in \cite[Section 5]{FeNo2}:
\begin{multline}\label{dispersive.1}
\sup_{t,\in[0,T]}\left( \|\nabla \Phi_{\varepsilon,\eta}(t,\cdot)\|_{W^{k,2}(\Omega)\cap W^{k,\infty}(\Omega)} + \left\|\left( \frac{\overline \theta}{\overline\varrho} R_{\varepsilon,\eta} + T_{\varepsilon,\eta} -F\right)(t,\cdot)\right\|_{W^{k,2}(\Omega)\cap W^{k,\infty}(\Omega)}\right)\\ \leq c(k,\eta)\left(\|\nabla \Phi_{0,\eta}\|_{L^2(\Omega)} + \left\|\frac{\overline\theta}{\overline\varrho} R_{0,\eta} + T_{0,\eta}-F\right\|_{L^2(\Omega)}\right),
\end{multline}
\begin{multline}\label{dispersive.2}
\int_0^T \left( \|\nabla \Phi_{\varepsilon,\eta}(t,\cdot)\|_{ W^{k,\infty}(\Omega)} + \left\|\left(\frac{\overline \theta}{\overline\varrho} R_{\varepsilon,\eta} + T_{\varepsilon,\eta}-F\right)(t,\cdot)\right\|_{ W^{k,\infty}(\Omega)}\right)\ {\rm d}t\\
\leq \omega(\varepsilon,\eta,k)\left(\|\nabla \Phi_{0,\eta}\|_{L^2(\Omega)} + \left\|\frac{\overline\theta}{\overline\varrho} R_{0,\eta} + T_{0,\eta}-F\right\|_{L^2(\Omega)}\right)
\end{multline}
with $\omega(\varepsilon,\eta,k)\to 0$ as $\varepsilon\to 0$ for any fixed $\eta>0$ and $k\geq 0$. Moreover, 
\begin{equation}\label{dispersive.3}
\sup_{t\in[0,T]} \left\|\frac{c_v\overline\varrho}{\overline\theta} T_{\varepsilon,\eta}(t,\cdot) - R_{\varepsilon,\eta}(t,\cdot)\right\|_{W^{k,q}(\Omega)}\leq c(\eta,k,F)\left(1+\left\|\frac{c_v\overline\varrho}{\overline\theta} T_{0,\eta} - R_{0,\eta}\right\|_{L^2(\Omega)}\right)
\end{equation}
for $k=0,1, \dots$ and $q \in [1,\infty]$.

Further, $R_{\varepsilon,\eta}\to R_\eta$ strongly in $L^\infty_{loc}((0,T];W^{k,p}(\Omega))$, $T_{\varepsilon,\eta}\to T_\eta$ strongly in $L^\infty_{loc}((0,T]; W^{k,p}(\Omega))$ for $p>2$ and $k\in \mathbb N$ and $T_\eta$ solves
\begin{equation}
(c_v+1)\left(\partial_t T_\eta + {\bf U}\cdot \nabla T_\eta\right) - {\bf U}\cdot \nabla F = 0
\end{equation}
with the initial data
\begin{equation}
T_\eta(0,\cdot) = \frac{\overline\theta}{c_v+1} \left( c_v \frac 1{\overline\theta} \left[\theta_0^{(1)}\right]_\eta - \frac 1{\overline\varrho} \left[\varrho_0^{(1)}\right]_\eta + \frac 1{\overline\theta}\left[ F\right]_\eta\right).
\end{equation}
Finally,  we observe that $T_\eta\to \Theta$ strongly in $L^\infty([0,T];L^2(\Omega))$  and  since $R_{\varepsilon,\eta}\to \frac{\overline\varrho}{\overline\theta} T_{\varepsilon,\eta} + \frac{\overline\varrho}{\overline\theta}F$ as $\varepsilon \to 0$ we conclude from \eqref{Euler.boussinesq}$_4$ that $R_{\eta}\to r$ 
strongly in $L^\infty([0,T];L^2(\Omega))$.

\subsection{Convergence in the relative energy inequality}

Due \eqref{rel.en.drop.chi}  and \eqref{rel.ini.ill} the relative energy inequality for \eqref{test.function.ill.posed} reads as follows
\begin{multline}
\int_\Omega \langle U^\varepsilon_{\tau,x};\mathcal E^\varepsilon(\varrho,{\bf m}, p|\tilde \varrho_{\varepsilon,\eta},\tilde{\bf U}_{\varepsilon,\eta},\tilde \theta_{\varepsilon,\eta})\rangle\ {\rm d}x + \mathcal D^\varepsilon(\tau)\\
\leq -\frac 1{\varepsilon}\int_0^\tau \int_\Omega \langle U^\varepsilon_{t,x};\varrho \chi_{a,b}(s(\varrho, p))\rangle \partial_t T_{\varepsilon,\eta} + \langle U_{t,x}^\varepsilon; \chi_{a,b}(s(\varrho, p)){\bf m}\rangle\cdot \nabla T_{\varepsilon,\eta}\ {\rm d}x{\rm d}t\\
+\int_0^\tau \int_\Omega \langle U_{t,x}^\varepsilon;\varrho \tilde{\bf U}_{\varepsilon,\eta} - {\bf m}\rangle \cdot \partial_t \tilde{\bf U}_{\varepsilon,\eta} + \left\langle U_{t,x}^\varepsilon;\frac{(\varrho\tilde{\bf U}_{\varepsilon,\eta} - {\bf m})\otimes{\bf m}}{\varrho}\right \rangle :\nabla \tilde{\bf U}_{\varepsilon,\eta} \ {\rm d}x{\rm d}t\\
- \frac 1{\varepsilon^2} \int_0^\tau \int_\Omega \langle U_{t,x}^\varepsilon; p\rangle \Delta \Phi_{\varepsilon,\eta}\ {\rm d}x{\rm d}t\\
+ \frac 1{\varepsilon} \int_0^\tau \int_\Omega \langle U_{t,x}^\varepsilon; \varrho \rangle \partial_t  T_{\varepsilon,\eta} s(\tilde \vr_{\varepsilon,\eta}, \tilde \varrho_{\varepsilon,\eta} \tilde \theta_{\varepsilon,\eta}) + \langle U_{t,x}^\varepsilon;{\bf m}\rangle \cdot \nabla  T_{\varepsilon,\eta} s(\tilde \vr_{\varepsilon,\eta},\tilde \varrho_{\varepsilon,\eta} \tilde \theta_{\varepsilon,\eta})\ {\rm d}x{\rm d}t\\
+\frac 1{\varepsilon^2}\int_0^\tau \int_\Omega \langle U_{t,x}^\varepsilon; \tilde\varrho_{\varepsilon,\eta}-\varrho\rangle \frac 1{\tilde \varrho_{\varepsilon,\eta}} \partial_t  (\tilde\varrho_{\varepsilon,\eta} \tilde \theta_{\varepsilon,\eta}) - \langle U_{t,x}^\varepsilon;{\bf m}\rangle \cdot \frac 1{\tilde\varrho_{\varepsilon,\eta}} \nabla (\tilde \varrho_{\varepsilon,\eta} \tilde \theta_{\varepsilon,\eta})\ {\rm d}x{\rm d}t\\
+\frac 1\varepsilon \int_0^\tau \int_\Omega  \langle U_{t,x}^\varepsilon;{\bf m} - \varrho \tilde{\bf U}_{\varepsilon,\eta}\rangle \cdot \nabla F\ {\rm d}x{\rm d}t + \int_0^\tau \int_\Omega \nabla \tilde{\bf U}_{\varepsilon,\eta}: {\rm d}\mu_{c,\varepsilon} + \omega(\varepsilon,\eta)\\
=: I + II + III + IV + V + VI + VII + \omega(\varepsilon,\eta)
\end{multline}
for a.a. $\tau \in (0,T)$.

The regularity properties of $\Phi_{\varepsilon,\eta}$ (see \cite[Section 5.2]{FeNo2}) implies
$$
VII \aleq \int_0^\tau \mathcal D^\varepsilon(t) \ {\rm d}t.
$$

Similarly as in the well-prepared case (see the comments under the relation \eqref{respes}) we control 
\begin{equation}
\frac1{\varepsilon} \int_0^\tau \int_\Omega \langle U_{t,x}^\varepsilon; ({\bf m} - \varrho \tilde {\bf U}_{\varepsilon, \eta}) (\chi_{a,b}(s(\varrho,p)) - \chi_{a,b} (s(\tilde\varrho_{\varepsilon,\eta},\tilde\varrho_{\varepsilon,\eta}\tilde\theta_{\varepsilon,\eta}))\rangle \cdot \nabla T_{\varepsilon,\eta}\ {\rm d}x{\rm d}t
\end{equation}
and thus  $ I + IV$ can be replaced by
\begin{equation}
-\frac 1\varepsilon\int_0^\tau \int_\Omega \langle U_{t,x}^\varepsilon; \varrho(\chi_{a,b}(s(\varrho,p)) - s(\tilde\varrho_{\varepsilon,\eta},\tilde\varrho_{\varepsilon,\eta}\tilde\theta_{\varepsilon,\eta})) \rangle\\ (\partial_t T_{\varepsilon,\eta} + \tilde{\bf U}_{\varepsilon,\eta} \cdot \nabla T_{\varepsilon,\eta})\ {\rm d}x{\rm d}t,
\end{equation}
$\omega(\varepsilon, \eta)$ and \eqref{entropy.on.the.right}. 

By \eqref{Euler.boussinesq}$_{(2)}$ we calculate
\begin{multline}
II = \int_0^\tau \int_{\Omega} \langle U_{t,x}^\varepsilon; \varrho \tilde{\bf U}_{\varepsilon,\eta} - {\bf m}\rangle \cdot  \partial_t \nabla \Phi_{\varepsilon,\eta} - \langle U_{t,x}^\varepsilon; \varrho \tilde{\bf U}_{\varepsilon,\eta} - {\bf m}\rangle \cdot \nabla \Pi + \langle U_{t,x}^\varepsilon; \varrho \tilde{\bf U}_{\varepsilon,\eta} - {\bf m}\rangle\cdot \frac{r}{\overline \varrho} \nabla F\\
 + \left \langle U_{t,x}^\varepsilon; \frac{(\varrho \tilde{\bf U}_{\varepsilon,\eta} - {\bf m})\otimes({\bf m} - \varrho \tilde{\bf U}_{\varepsilon,\eta})}{\varrho}\right\rangle :\nabla \tilde{\bf U}_{\varepsilon,\eta} + \langle U_{t,x}^\varepsilon; (\varrho \tilde{\bf U}_{\varepsilon,\eta} - {\bf m})\otimes \nabla \Phi_{\varepsilon,\eta} \rangle : \nabla \tilde{\bf U}_{\varepsilon,\eta}\\
 + \langle U_{t,x}^\varepsilon;( \varrho\tilde{\bf U}_{\varepsilon,\eta} - {\bf m})\otimes {\bf U}_{\varepsilon,\eta}\rangle:\nabla^2\Phi_{\varepsilon,\eta} \ {\rm d}x{\rm d}t.
\end{multline}
The second and fourth terms can be included into $\omega(\varepsilon,\eta)$ and
\begin{equation}\label{entropy.on.the.right}
\int_0^\tau \int_\Omega \langle U^\varepsilon_{t,x};\mathcal E^\varepsilon(\varrho,{\bf m}, p|\tilde \varrho_{\varepsilon,\eta},\tilde{\bf U}_{\varepsilon,\eta},\tilde \theta_{\varepsilon,\eta})\rangle\ {\rm d}x{\rm d}t
\end{equation} 
due to \eqref{dispersive.1} and \eqref{dispersive.2}.  This can be done analogously to the well-prepared case.
The fifth and sixth term are handled similarly by the Young inequality and  \eqref{dispersive.1} to get
\begin{multline}
\int_0^\tau \int_\Omega\langle U_{t,x}^\varepsilon; (\varrho \tilde{\bf U}_{\varepsilon,\eta} - {\bf m})\otimes \nabla \Phi_{\varepsilon,\eta} \rangle : \nabla \tilde{\bf U}_{\varepsilon,\eta}\ {\rm d}x{\rm d}t \\
\aleq \int_0^\tau \int_\Omega \left\langle U_{t,x}^\varepsilon; \varrho \left|\frac{\bf m}\varrho - \tilde{\bf U}_{\varepsilon,\eta}\right|^2\right\rangle \ {\rm d}x{\rm d}t + \int_0^\tau \int_\Omega \langle U_{t,x}^\varepsilon; \varrho\rangle |\nabla \Phi_{\varepsilon, \eta}|^2\ {\rm d}x{\rm d}t\\
\aleq \int_0^\tau \int_\Omega \langle U_{t,x}^\varepsilon; \mathcal E^\varepsilon (\varrho,{\bf m},p|\tilde \varrho_{\varepsilon,\eta}, \tilde{\bf U}_{\varepsilon,\eta}, \tilde{\theta}_{\varepsilon,\eta})\rangle\  {\rm d}x{\rm d}t + \omega(\varepsilon,\eta).
\end{multline}
Finally, due to \eqref{acoustic.equation}$_2$, \eqref{rel.est.2} and the dispersive estimate \eqref{dispersive.2}
\begin{multline}
\int_0^\tau \int_{\Omega} \langle U_{t,x}^\varepsilon; \varrho \tilde{\bf U}_{\varepsilon,\eta} - {\bf m}\rangle \cdot \partial_t \nabla \Phi_{\varepsilon,\eta}\ {\rm d} x{\rm d}t \\
= \int_0^\tau \int_\Omega \langle U_{t,x}^\varepsilon; \varrho {\bf U} - {\bf m}\rangle \cdot \partial_t \nabla \Phi_{\varepsilon,\eta}\ {\rm d}x{\rm d}t + \int_0^\tau\int_\Omega \langle U_{t,x}^\varepsilon; \varrho \nabla \Phi_{\varepsilon,\eta}\rangle \cdot  \partial_t \nabla \Phi_{\varepsilon,\eta}\ {\rm d}x{\rm d}t\\
 = -\int_0^\tau \int_\Omega \langle U_{t,x}^\varepsilon; {\bf m} \rangle \cdot  \partial_t \nabla\Phi_{\varepsilon,\eta}\ {\rm d}x{\rm d}t + \frac 12 \int_0^\tau \int_\Omega \overline\varrho \partial_t |\nabla\Phi_{\varepsilon,\eta}|^2\ {\rm d}x {\rm d}t\\
-\int_0^\tau \int_\Omega\left \langle U_{t,x}^\varepsilon; \frac{\varrho - \overline\varrho}{\varepsilon}\right\rangle \nabla \Phi_{\varepsilon,\eta} \cdot  \nabla\left(\frac{\overline\theta}{\overline\varrho} R_{\varepsilon,\eta} + T_{\varepsilon,\eta} - F\right)\ {\rm d}x{\rm d}t\\
- \int_0^\tau \int_\Omega\left\langle U_{t,x}^\varepsilon; \frac{\varrho - \overline\varrho}{\varepsilon}  \right\rangle {\bf U} \cdot\nabla\left(\frac{\overline\theta}{\overline\varrho}R_{\varepsilon,\eta} + T_{\varepsilon,\eta} - F\right) \ {\rm d}x{\rm d}t \\
= -\int_0^\tau \int_\Omega \langle U_{t,x}^\varepsilon; {\bf m}\rangle\cdot  \partial_t \nabla \Phi_{\varepsilon,\eta} \ {\rm d}x{\rm d}t + \frac 12 \int_0^\tau \int_\Omega \overline\varrho \partial_t |\nabla\Phi_{\varepsilon,\eta}|^2\ {\rm d}x {\rm d}t  + \omega(\varepsilon,\eta).
\end{multline}

Put all together, we get
\begin{multline}\label{term.ii}
II = VIII  + \omega(\varepsilon,\eta) -\int_0^\tau \int_\Omega \langle U_{t,x}^\varepsilon; {\bf m}\rangle \cdot  \partial_t \nabla \Phi_{\varepsilon,\eta} \ {\rm d}x{\rm d}t \\
+ \int_0^\tau \int_\Omega \langle U_{t,x}^\varepsilon; \varrho \tilde{\bf U}_{\varepsilon,\eta} - {\bf m}\rangle\cdot \frac{r}{\overline \varrho} \nabla F \ {\rm d}x{\rm d}t + \frac 12 \int_0^\tau \int_\Omega \overline\varrho \partial_t |\nabla\Phi_{\varepsilon,\eta}|^2\ {\rm d}x {\rm d}t,
\end{multline}
where \begin{equation}VIII\aleq \int_0^\tau \int_\Omega\langle U_{t,x}^\varepsilon; \mathcal E^\varepsilon(\varrho,{\bf m},p|\tilde{\varrho}_{\varepsilon,\eta}, \tilde{\bf U}_{\varepsilon,\eta},\tilde{\theta}_{\varepsilon,\eta})\rangle \ {\rm d}x{\rm d}t.\end{equation}

We handle the second term in $V$  using \eqref{acoustic.equation}$_{(2)}$ and the dispersive estimate \eqref{dispersive.2} as
\begin{multline}
-\frac 1{\varepsilon^2} \int_0^\tau \int_\Omega \langle U_{t,x}^\varepsilon; {\bf m}\rangle \cdot\frac 1{\tilde\varrho_{\varepsilon,\eta}}\nabla(\tilde\varrho_{\varepsilon,\eta} \tilde \theta_{\varepsilon,\eta})\ {\rm d}x{\rm d}t \\
=  - \frac 1\varepsilon\int_0^\tau \int_\Omega \langle U_{t,x}^\varepsilon,{\bf m}\rangle  \cdot\frac 1{\tilde\varrho_{\varepsilon,\eta}} \nabla(\overline\varrho T_{\varepsilon,\eta} + \overline\theta R_{\varepsilon,\eta})\ {\rm d}x{\rm d}t\\
 -  \int_0^\tau \int_\Omega\langle U_{t,x}^\varepsilon;{\bf m}\rangle\cdot \frac 1{\tilde\varrho_{\varepsilon,\eta}}\nabla (R_{\varepsilon,\eta}T_{\varepsilon,\eta})\ {\rm d}x{\rm d}t \\
=  - \frac 1\varepsilon\int_0^\tau \int_\Omega\langle U_{t,x}^\varepsilon;{\bf m} \rangle\cdot \frac{\overline\varrho}{\tilde\varrho_{\varepsilon,\eta}}\nabla\left(T_{\varepsilon,\eta} + \frac{\overline\theta}{\overline\varrho} R_{\varepsilon,\eta}\right)\ {\rm d}x{\rm d}t + \omega(\varepsilon,\eta)\\
 =-\frac 1{\varepsilon}\int_0^\tau \int_\Omega \langle U_{t,x}^\varepsilon,{\bf m}\rangle\cdot \nabla F\ {\rm d}x{\rm d}t + \int_0^\tau \int_\Omega \langle U_{t,x}^\varepsilon,{\bf m}\rangle\cdot \frac {R_{\varepsilon,\eta}}{\tilde \varrho_{\varepsilon,\eta}} \nabla F\ {\rm d}x{\rm d}t\\
 + \int_0^\tau\int_\Omega\langle U_{t,x}^\varepsilon;{\bf m}\rangle\cdot \partial_t \nabla \Phi_{\varepsilon, \eta} \ {\rm d}x{\rm d}t + \omega(\varepsilon,\eta),
\end{multline}
where the third term on the right-hand side cancels with a term from $II$ (see \eqref{term.ii}) and we used \eqref{solenoidal.limit} to calculate
\begin{multline}
\int_0^\tau \int_\Omega \langle U_{t,x}^\varepsilon,{\bf m}\rangle\cdot \frac 1{\overline\varrho}\frac {\overline\varrho}{\tilde \varrho_{\varepsilon,\eta}} \nabla\left(R_{\varepsilon,\eta}T_{\varepsilon,\eta}\right)\ {\rm d}x{\rm d}t\\
= \int_0^\tau \int_\Omega \langle U_{t,x}^\varepsilon; {\bf m}\rangle \cdot\frac 1{\overline\varrho} \nabla \left(R_{\varepsilon, \eta} T_{\varepsilon,\eta}\right)\ {\rm d}x{\rm d}t - \int_0^\tau \int_\Omega \langle U_{t,x}^\varepsilon ; {\bf m}\rangle\cdot \frac 1{\overline\varrho}\frac{\varepsilon R_{\varepsilon,\eta}}{\tilde \varrho_{\varepsilon,\eta}} \nabla \left(R_{\varepsilon,\eta} T_{\varepsilon,\eta}\right)\ {\rm d}x{\rm d}t\\
=\omega(\varepsilon,\eta).
\end{multline}
Concerning the first term in $V$, we proceed as follows
\begin{multline}
\frac 1{\varepsilon^2}\int_0^\tau \int_\Omega \langle U_{t,x}^\varepsilon; \tilde\varrho_{\varepsilon,\eta}-\varrho\rangle \frac 1{{\tilde \varrho}_{\varepsilon, \eta}} \partial_t  (\tilde\varrho_{\varepsilon,\eta} \tilde \theta_{\varepsilon,\eta})\ {\rm d}x{\rm d}t\\
 = \int_0^\tau \int_\Omega \left\langle U_{t,x}^\varepsilon; \frac{\tilde \varrho_{\varepsilon,\eta} - \varrho}{\varepsilon\tilde \varrho_{\varepsilon,\eta}}\right\rangle \left(\overline\varrho \partial_t T_{\varepsilon,\eta} + \overline \theta \partial_t R_{\varepsilon,\eta}\right)\ {\rm d}x{\rm d}t\\
+ \int_0^\tau \int_\Omega \left\langle U_{t,x}^\varepsilon; \frac{\tilde\varrho_{\varepsilon,\eta} - \varrho}{\varepsilon \tilde \varrho_{\varepsilon,\eta}}\right\rangle \varepsilon\partial_t\left(R_{\varepsilon,\eta} T_{\varepsilon,\eta}\right)\ {\rm d}x {\rm d}t,
\end{multline}
where the last term on the right-hand side can be included into $\omega(\varepsilon,\eta)$. Indeed, we may deduce from \eqref{acoustic.equation} and \eqref{temperature.for.T} that
\begin{equation}\label{time.derivative.T}
\varepsilon\partial_t T_{\varepsilon,\eta} = -\frac{\overline\theta}{c_v} \Delta \Phi_{\varepsilon,\eta} - \frac{\varepsilon}{c_v+1} \tilde{\bf U}_{\varepsilon,\eta}\cdot \nabla \left( c_v T_\varepsilon - \frac{\overline\theta}{\overline\varrho} R_\varepsilon\right)
\end{equation}
and the rest follows from \eqref{dispersive.1} and \eqref{dispersive.2}. Similar computation may be deduced even for $R_{\varepsilon,\eta}$.

Since $R_{\varepsilon, \eta} \to r$ in $L^\infty_{loc}((0,T];L^2(\Omega))$ we obtain
\begin{equation}\label{boussinesq.limit}
\int_0^\tau \int_\Omega \langle U_{t,x}^\varepsilon;{\bf m}\rangle \cdot \frac{R_{\varepsilon,\eta}}{\tilde \varrho_{\varepsilon,\eta}} \nabla F\ {\rm d}x{\rm d}t = \int_0^\tau \int_\Omega \langle U_{t,x}^\varepsilon;{\bf m} \rangle\cdot \frac r{\overline\varrho} \nabla F\ {\rm d}x{\rm d}t + \omega(\varepsilon,\eta) + \mathcal E.
\end{equation}
where
\begin{equation}
\mathcal E\aleq \int_0^\tau \int_{\Omega} \langle U_{t,x}^\varepsilon; \mathcal E^\varepsilon (\varrho,{\bf m},p|\tilde \varrho_{\varepsilon,\eta}, \tilde{\bf U}_{\varepsilon,\eta}, \tilde \theta_{\varepsilon,\eta}\rangle \ {\rm d}x {\rm d}t.
\end{equation}
We collect all the computations in order to get
\begin{multline}\label{mezikrok}
\int_\Omega \langle U^\varepsilon_{\tau,x};\mathcal E^\varepsilon(\varrho,{\bf m}, p|\tilde \varrho_{\varepsilon,\eta},\tilde{\bf U}_{\varepsilon,\eta},\tilde \theta_{\varepsilon,\eta})\rangle\ {\rm d}x + \mathcal D^\varepsilon(\tau)\\
\aleq \int_0^\tau \int_\Omega \langle U^\varepsilon_{t,x};\mathcal E^\varepsilon(\varrho,{\bf m}, p|\tilde \varrho_{\varepsilon,\eta},\tilde{\bf U}_{\varepsilon,\eta},\tilde \theta_{\varepsilon,\eta})\rangle\ + \int_0^\tau \mathcal D^\varepsilon(t)\ {\rm d}t+\omega(\varepsilon,\eta)\\
 - \frac 1\varepsilon\int_0^\tau \int_\Omega \langle U_{t,x}^\varepsilon; \varrho(\chi(s(\varrho,p)) - s(\tilde\varrho_{\varepsilon,\eta},\tilde\varrho_{\varepsilon,\eta}\tilde\theta_{\varepsilon,\eta})) \rangle  (\partial_t T_{\varepsilon,\eta} + \tilde{\bf U}_{\varepsilon,\eta} \cdot\nabla T_{\varepsilon,\eta})\ {\rm d}x{\rm d}t\\
+\int_0^\tau \int_\Omega \langle U_{t,x}^\varepsilon; \varrho\tilde{\bf U}_{\varepsilon,\eta}\rangle \cdot \frac r{\overline\varrho} \nabla F\ {\rm d}x{\rm d}t
- \frac 1{\varepsilon^2} \int_0^\tau \int_\Omega \langle U_{t,x}^\varepsilon; p\rangle \Delta \Phi_{\varepsilon,\eta}\ {\rm d}x{\rm d}t\\
+\int_0^\tau \int_\Omega \left\langle U_{t,x}^\varepsilon; \frac{\tilde \varrho_{\varepsilon,\eta} - \varrho}{\varepsilon\tilde \varrho_{\varepsilon,\eta}}\right\rangle \left(\overline\varrho \partial_t T_{\varepsilon,\eta} + \overline \theta \partial_t R_{\varepsilon,\eta}\right)\ {\rm d}x{\rm d}t\\
+ \frac 12 \int_0^\tau \int_\Omega \overline\varrho \partial_t |\nabla\Phi_{\varepsilon,\eta}|^2\ {\rm d}x {\rm d}t -\frac 1\varepsilon \int_0^\tau \int_\Omega  \langle U_{t,x}^\varepsilon;\varrho \tilde{\bf U}_{\varepsilon, \eta}\rangle \cdot \nabla F\ {\rm d}x{\rm d}t
\end{multline}
for a.a. $\tau \in (0,T)$.

The residual part of the third term on the right-hand side is handled with help of \eqref{rel.est}, \eqref{dispersive.1} and \eqref{dispersive.2} as follows:
\begin{multline}
-\frac 1\varepsilon\int_0^\tau \int_\Omega \langle U_{t,x}^\varepsilon;\varrho_{res} (\chi(s(\varrho,p)) - s(\tilde\varrho_{\varepsilon,\eta},\tilde\varrho_{\varepsilon,\eta}\tilde\theta_{\varepsilon,\eta}))_{res} \rangle (\partial_t T_{\varepsilon,\eta} + \tilde{\bf U}_{\varepsilon,\eta} \cdot \nabla T_{\varepsilon,\eta})\ {\rm d}x{\rm d}t\\
\aleq \left|\int_0^\tau \int_\Omega \left\langle U_{t,x}^\varepsilon; \frac {\varrho_{res}}{\varepsilon^2}\right\rangle \varepsilon \left(\partial_t T_{\varepsilon,\eta} + \tilde{\bf U}_{\varepsilon,\eta}\cdot \nabla T_{\varepsilon,\eta}\right)\ {\rm d}x{\rm d}t\right|\\
\aleq \int_0^\tau \int_\Omega \langle U_{t,x}^\varepsilon; \mathcal E^\varepsilon (\varrho, {\bf m}, p|\tilde{\varrho}_{\varepsilon,\eta}, \tilde{\bf U}_{\varepsilon,\eta}, \tilde\theta_{\varepsilon,\eta}\rangle\ {\rm d}x{\rm d}t.
\end{multline}
We use the Taylor formula to manage the essential part similarly as in \eqref{taylor} to get
\begin{multline}
-\frac 1\varepsilon\int_0^\tau \int_\Omega \langle U_{t,x}^\varepsilon;\varrho_{ess} (\chi(s(\varrho,p)) - s(\tilde\varrho_{\varepsilon,\eta},\tilde\varrho_{\varepsilon,\eta}\tilde\theta_{\varepsilon,\eta}))_{ess} \rangle (\partial_t T_{\varepsilon,\eta} + \tilde{\bf U}_{\varepsilon,\eta}\cdot \nabla T_{\varepsilon,\eta})\ {\rm d}x{\rm d}t\\
= \frac{c_v + 1}{\varepsilon}\int_0^\tau \int_\Omega \left\langle U_{t,x}^\varepsilon; \left(\frac{\varrho - \tilde \varrho_{\varepsilon,\eta}}{\tilde\varrho_{\varepsilon,\eta}} 	 \varrho\right)_{ess}\right\rangle \left(\partial_t T_{\varepsilon,\eta} + \tilde{\bf U}_{\varepsilon,\eta}\cdot\nabla T_{\varepsilon,\eta}\right)\ {\rm d}x{\rm d}t\\
-\frac{c_v}{\varepsilon} \int_0^\tau \int_\Omega \left\langle U_{t,x}^\varepsilon; \left(\frac{p - \tilde\varrho_{\varepsilon,\eta}\tilde\theta_{\varepsilon,\eta}}{\tilde\varrho_{\varepsilon,\eta}\tilde\theta_{\varepsilon,\eta}}\varrho\right)_{ess}\right\rangle \left(\partial_t T_{\varepsilon,\eta} + \tilde{\bf U}_{\varepsilon,\eta} \cdot\nabla T_{\varepsilon,\eta}\right)\ {\rm d}x{\rm d}t + \omega(\varepsilon,\eta)\\
= \frac{c_v + 1}{\varepsilon}\int_0^\tau \int_\Omega \left\langle U_{t,x}^\varepsilon; \frac{\varrho - \tilde \varrho_{\varepsilon,\eta}}{\tilde\varrho_{\varepsilon,\eta}} 	 \overline\varrho\right\rangle \left(\partial_t T_{\varepsilon,\eta} + \tilde{\bf U}_{\varepsilon,\eta}\cdot\nabla T_{\varepsilon,\eta}\right)\ {\rm d}x{\rm d}t\\
-\frac{c_v}{\varepsilon} \int_0^\tau \int_\Omega \left\langle U_{t,x}^\varepsilon; \frac{p - \tilde\varrho_{\varepsilon,\eta}\tilde\theta_{\varepsilon,\eta}}{\overline\theta}\right\rangle \left(\partial_t T_{\varepsilon,\eta} + \tilde{\bf U}_{\varepsilon,\eta} \nabla T_{\varepsilon,\eta}\right)\ {\rm d}x{\rm d}t + \omega(\varepsilon,\eta)
\end{multline} 
We use \eqref{acoustic.equation}$_1$ to get
\begin{multline}
-\frac 1{\varepsilon^2}\int_0^\tau \int_\Omega \langle U_{t,x}^\varepsilon; p\rangle \Delta \Phi_{\varepsilon,\eta}\ {\rm d}x{\rm d}t = -\frac 1{\varepsilon^2}\int_0^\tau \int_\Omega \langle U_{t,x}^\varepsilon; p - \overline\varrho\overline\theta\rangle \Delta \Phi_{\varepsilon,\eta}\ {\rm d}x{\rm d}t\\
= \frac 1{\varepsilon}\int_0^\tau \int_\Omega \langle U_{t,x}^\varepsilon; p-\overline\varrho \overline\theta\rangle \frac{c_v}{c_v+1} \partial_t \left( \frac {R_{\varepsilon,\eta}}{\overline\varrho} + \frac{T_{\varepsilon,\eta}}{\overline\theta}\right)\ {\rm d}x{\rm d}t.
\end{multline}
According to \eqref{test.function.ill.posed}
\begin{multline}
\int_0^\tau \int_\Omega \left\langle U_{t,x}^\varepsilon, \frac{p - \overline\varrho \overline\theta}{\varepsilon}\right\rangle \frac{c_v}{c_v + 1}\partial_t \left(\frac{R_{\varepsilon,\eta}}{\overline\varrho} + \frac{T_{\varepsilon,\eta}}{\overline\theta}\right)\ {\rm d}x{\rm d}t \\
 = \int_0^\tau \int_\Omega \left\langle U_{t,x}^\varepsilon, \frac{p - \tilde \varrho_{\varepsilon,\eta}\tilde \theta_{\varepsilon,\eta}}{\varepsilon}\right\rangle\frac{c_v}{c_v+1} \partial_t  \left( \frac{R_{\varepsilon,\eta}}{\overline\varrho} + \frac{T_{\varepsilon,\eta}}{\overline\theta}\right)\ {\rm d}x{\rm d}t\\
+\int_0^\tau \int_\Omega \left(\overline\theta R_{\varepsilon,\eta} + \overline\varrho T_{\varepsilon,\eta}\right) \frac{c_v}{c_v+1}\partial_t\left(\frac{R_{\varepsilon,\eta}}{\overline\varrho} + \frac{T_{\varepsilon,\eta}}{\overline\theta}\right)\ {\rm d}x{\rm d}t\\
+ \int_0^\tau \int_{\Omega} R_{\varepsilon,\eta}T_{\varepsilon,\eta} \frac{c_v}{c_v+1} \varepsilon\partial_t\left(\frac{R_{\varepsilon,\eta}}{\overline\varrho} + \frac{T_{\varepsilon,\eta}}{\overline\theta}\right)\ {\rm d}x{\rm d}t,
\end{multline}
where the last term can be included into $\omega(\varepsilon,\eta)$ due to \eqref{time.derivative.T}, \eqref{dispersive.1} and \eqref{dispersive.2}.

Consequently, we get from \eqref{mezikrok} that
\begin{multline}\label{mezikrok.2}
\int_\Omega \langle U^\varepsilon_{\tau,x};\mathcal E^\varepsilon(\varrho,{\bf m}, p|\tilde \varrho_{\varepsilon,\eta},\tilde{\bf U}_{\varepsilon,\eta},\tilde \theta_{\varepsilon,\eta})\rangle\ {\rm d}x + \mathcal D^\varepsilon(\tau)\\
\aleq \int_0^\tau \int_\Omega \langle U^\varepsilon_{t,x};\mathcal E^\varepsilon(\varrho,{\bf m}, p|\tilde \varrho_{\varepsilon,\eta},\tilde{\bf U}_{\varepsilon,\eta},\tilde \theta_{\varepsilon,\eta})\rangle\ + \int_0^\tau \mathcal D^\varepsilon(t)\ {\rm d}t+\omega(\varepsilon,\eta)\\
+ 
\int_0^\tau \int_\Omega \langle U_{t,x}^\varepsilon; \varrho\tilde{\bf U}_{\varepsilon,\eta}\rangle \cdot  \frac r{\overline\varrho} \nabla F\ {\rm d}x{\rm d}t-\frac 1\varepsilon \int_0^\tau \int_\Omega  \langle U_{t,x}^\varepsilon; \varrho \tilde{\bf U}_{\varepsilon, \eta}\rangle \cdot \nabla F\ {\rm d}x{\rm d}t+ \frac 12 \int_0^\tau \int_\Omega \overline\varrho \partial_t |\nabla\Phi_{\varepsilon,\eta}|^2\ {\rm d}x {\rm d}t\\
+\int_0^\tau \int_\Omega \left\langle U_{t,x}^\varepsilon; \frac{\tilde \varrho_{\varepsilon,\eta} - \varrho}{\varepsilon\tilde \varrho_{\varepsilon,\eta}}\right\rangle \left(\left(\overline\varrho\partial_t T_{\varepsilon,\eta} + \overline \theta \partial_t R_{\varepsilon,\eta}\right) - (c_v + 1)\left(\overline\varrho\partial_t T_{\varepsilon,\eta} + \overline\varrho \tilde {\bf U}_{\varepsilon,\eta} \cdot\nabla T_{\varepsilon,\eta}\right)\right)\ {\rm d}x{\rm d}t\\
 + \int_0^\tau \int_{\Omega} \left\langle U_{t,x}^\varepsilon; \frac{p-\tilde\varrho_{\varepsilon,\eta}\tilde\theta_{\varepsilon,\eta}}{\varepsilon\overline\theta}\right\rangle \left(-c_v\left(\partial_t T_{\varepsilon,\eta} + \tilde{\bf U}_{\varepsilon,\eta}\cdot \nabla T_{\varepsilon,\eta}\right) + \frac{c_v}{c_v +1} \partial_t \left(\frac{\overline\theta}{\overline\rho} R_{\varepsilon,\eta} + T_{\varepsilon,\eta}\right)\right)\ {\rm d}x{\rm d}t\\
+\int_0^\tau \int_{\Omega}\frac{c_v}{c_v+1} \left(\overline\theta R_{\varepsilon,\eta} + \overline\varrho T_{\varepsilon,\eta}\right)\partial_t\left(\frac{R_{\varepsilon,\eta}}{\overline\varrho} + \frac{T_{\varepsilon,\eta}}{\overline\theta}\right)\ {\rm d}x{\rm d}t
\end{multline}
for a.a. $\tau \in (0,T)$.

We compute
\begin{multline}
\overline\varrho\partial_t T_{\varepsilon,\eta} + \overline\theta \partial_t R_{\varepsilon,\eta} - (c_v+1) \left(\overline\varrho\partial_t T_{\varepsilon,\eta}  + \overline\varrho \tilde{\bf U}_{\varepsilon,\eta}\cdot \nabla T_{\varepsilon,\eta}\right) \\
= -\left( c_v\overline\varrho \partial_t T_{\varepsilon,\eta} - \overline\theta \partial_t R_{\varepsilon,\eta} + \tilde{\bf U}_{\varepsilon,\eta} \cdot (c_v\overline\varrho\nabla T_{\varepsilon,\eta} - \overline\theta\nabla R_{\varepsilon,\eta})\right)\\
 - \tilde{\bf U}_{\varepsilon,\eta}\cdot  (\overline\varrho \nabla T_{\varepsilon,\eta} + \overline \theta \nabla R_{\varepsilon,\eta}-\overline\varrho\nabla F) - \tilde{\bf U}_{\varepsilon,\eta}\cdot\overline \varrho\nabla F,
\end{multline}
where the first term equals $0$ due to \eqref{temperature.for.T} and the second term tends to $0$ in the appropriate integrals due to the dispersive estimate \eqref{dispersive.2}.

Thus we rewrite \eqref{mezikrok.2} as
\begin{multline}\label{mezikrok.3}
\int_\Omega \langle U^\varepsilon_{\tau,x};\mathcal E^\varepsilon(\varrho,{\bf m}, p|\tilde \varrho_{\varepsilon,\eta},\tilde{\bf U}_{\varepsilon,\eta},\tilde \theta_{\varepsilon,\eta})\rangle\ {\rm d}x + \mathcal D^\varepsilon(\tau)\\
\leq \int_0^\tau \int_\Omega \langle U_{t,x}^\varepsilon; \varrho\tilde{\bf U}_{\varepsilon,\eta}\rangle\cdot \frac r{\overline\varrho} \nabla F\ {\rm d}x{\rm d}t-\frac 1\varepsilon \int_0^\tau \int_\Omega  \langle U_{t,x}^\varepsilon; \varrho \tilde{\bf U}_{\varepsilon,\eta}\rangle \cdot \nabla F\ {\rm d}x{\rm d}t\\
-\int_0^\tau \int_\Omega \left\langle U_{t,x}^\varepsilon; \frac{\tilde \varrho_{\varepsilon,\eta} - \varrho}{\varepsilon\tilde \varrho_{\varepsilon,\eta}}\right\rangle \tilde{\bf U}_{\varepsilon,\eta} \cdot\overline\varrho \nabla F\ {\rm d}x{\rm d}t+ \frac 12 \int_0^\tau \int_\Omega \overline\varrho \partial_t |\nabla\Phi_{\varepsilon,\eta}|^2\ {\rm d}x {\rm d}t\\
 - \int_0^\tau \int_{\Omega} \left\langle U_{t,x}^\varepsilon; \frac{p-\tilde\varrho_{\varepsilon,\eta}\tilde\theta_{\varepsilon,\eta}}{\varepsilon\overline\theta}\right\rangle \frac{c_v}{c_v+1} \tilde{\bf U}_{\varepsilon,\eta}\cdot\nabla F\ {\rm d}x{\rm d}t\\
+\int_0^\tau \int_{\Omega}\frac{c_v}{c_v+1} \left(\overline\theta R_{\varepsilon,\eta} + \overline\varrho T_{\varepsilon,\eta}\right)\partial_t\left(\frac{R_{\varepsilon,\eta}}{\overline\varrho} + \frac{T_{\varepsilon,\eta}}{\overline\theta}\right)\ {\rm d}x{\rm d}t\\
+\int_0^\tau \int_\Omega \langle U^\varepsilon_{t,x};\mathcal E^\varepsilon(\varrho,{\bf m}, p|\tilde \varrho_{\varepsilon,\eta},\tilde{\bf U}_{\varepsilon,\eta},\tilde \theta_{\varepsilon,\eta})\rangle\ + \int_0^\tau \mathcal D^\varepsilon(t)\ {\rm d}t+\omega(\varepsilon,\eta)
\end{multline}
for a.a. $\tau \in (0,T)$.

As
\begin{equation}
\tilde \varrho_{\varepsilon,\eta} - \varrho = \varepsilon R_{\varepsilon,\eta} + \overline \varrho - \varrho,
\end{equation}
the third term on the right-hand side is handled as follows
\begin{multline}
- \int_0^\tau\int_{\Omega} \left \langle U_{t,x}^\varepsilon; \frac{\tilde\varrho_{\varepsilon,\eta} - \varrho}{\varepsilon \tilde \varrho_{\varepsilon,\eta}}\right \rangle \tilde{\bf U}_{\varepsilon,\eta}\cdot \overline \varrho \nabla F\ {\rm d}x{\rm d}t \\
=-\int_0^\tau \int_\Omega  \frac1\varepsilon  \tilde{\bf U}_{\varepsilon,\eta}\cdot  \overline\varrho \nabla F\ {\rm d}x{\rm d}t + \int_0^\tau \int_\Omega \left\langle U_{t,x}^\varepsilon; \varrho \right \rangle \frac{\overline\varrho}{\varepsilon\tilde \varrho_{\varepsilon,\eta}} \tilde{\bf U}_{\varepsilon,\eta} \cdot \nabla F\ {\rm d}x{\rm d}t\\
=-\int_0^\tau \int_\Omega  \frac1\varepsilon  \left( U + \nabla \Phi_{\varepsilon,\eta}\right)\cdot  \overline\varrho \nabla F\ {\rm d}x{\rm d}t
  + \int_0^\tau \int_\Omega \langle U_{t,x}^\varepsilon; \varrho\rangle \frac 1\varepsilon \tilde{\bf U}_{\varepsilon,\eta} \cdot \nabla F\ {\rm d}x{\rm d}t\\
-  \int_0^\tau \int_\Omega \langle U_{t,x}^\varepsilon;\varrho\rangle \frac{R_{\varepsilon,\eta}}{\tilde\varrho_{\varepsilon,\eta}} \tilde{\bf U}_{\varepsilon,\eta}\cdot  \nabla F\ {\rm d}x{\rm d}t \\ 
=\int_0^\tau \int_\Omega  \frac1\varepsilon  \Delta \Phi_{\varepsilon,\eta}\cdot  \overline\varrho F\ {\rm d}x{\rm d}t
 + \int_0^\tau \int_\Omega \langle U_{t,x}^\varepsilon; \varrho\rangle \frac 1\varepsilon \tilde{\bf U}_{\varepsilon,\eta} \cdot \nabla F\ {\rm d}x{\rm d}t\\
- \int_0^\tau \int_\Omega \langle U_{t,x}^\varepsilon, \varrho\rangle \frac r{\overline\varrho} \tilde{\bf U}_{\varepsilon,\eta}\cdot \nabla F\ {\rm d}x{\rm d}t  + \omega(\varepsilon, \eta)\\
= - \int_0^\tau \int_\Omega \langle U_{t,x}^\varepsilon, \varrho\rangle \frac r{\overline\varrho} \tilde{\bf U}_{\varepsilon,\eta}\cdot \nabla F\ {\rm d}x{\rm d}t + \int_0^\tau \int_\Omega \langle U_{t,x}^\varepsilon; \varrho\rangle \frac 1\varepsilon \tilde{\bf U}_{\varepsilon,\eta} \cdot \nabla F\ {\rm d}x{\rm d}t\\
- \int_0^\tau \int_{\Omega} \frac{c_v}{c_v+1} \partial_t \left(\frac{ R_{\varepsilon,\eta}}{\overline\varrho} + \frac{T_{\varepsilon,\eta}}{\overline\theta} - \frac{F}{\overline\theta}\right)\overline \varrho F\ {\rm d}x{\rm d}t + \omega(\varepsilon,\eta),
\end{multline}
where the first and the second term on the right-hand side cancels with the first two terms on the right-hand side of \eqref{mezikrok.3}.
Further, we compute
\begin{multline}
\int_0^\tau \int_{\Omega} \left \langle U_{t,x}^\varepsilon; \frac{p- \tilde\varrho_{\varepsilon,\eta}\tilde\theta_{\varepsilon,\eta}}{\varepsilon\overline\theta}\right\rangle \tilde{\bf U}_{\varepsilon,\eta}\cdot  \nabla F\ {\rm d}x{\rm d}t = \int_0^\tau\int_{\Omega} \left\langle U_{t,x}^\varepsilon; \frac{p - \tilde{\varrho}_{\varepsilon,\eta}\tilde\theta_{\varepsilon,\eta}}{\varepsilon\overline\theta}\right\rangle {\bf U} \cdot \nabla F \ {\rm d}x{\rm d}t + \omega(\varepsilon,\eta)\\
 =  \int_0^\tau \int_{\Omega} \left\langle U_{t,x}^\varepsilon; \frac{p-\overline\varrho\overline\theta}{\varepsilon\overline\theta}\right\rangle {\bf U}\cdot \nabla F\ {\rm d}x{\rm d}t -\int_0^\tau \int_\Omega \left(\overline\varrho T_{\varepsilon,\eta} + \overline\theta R_{\varepsilon,\eta}\right){\bf U}\cdot \nabla F \ {\rm d}x{\rm d}t + \omega(\varepsilon,\eta)\\
=  \int_0^\tau \int_{\Omega} \left\langle U_{t,x}^\varepsilon; \frac{p-\overline\varrho\overline\theta}{\varepsilon\overline\theta}\right\rangle {\bf U}\cdot \nabla F\ {\rm d}x{\rm d}t -\int_0^\tau \int_\Omega \left(\overline\varrho T_{\varepsilon,\eta} + \overline\theta R_{\varepsilon,\eta} - \overline\varrho F\right){\bf U}\cdot \nabla F \ {\rm d}x{\rm d}t  \\
+ \omega(\varepsilon,\eta)
=\int_0^\tau \int_{\Omega} \left\langle U_{t,x}^\varepsilon; \frac{p- \varepsilon \overline \varrho \Theta - \varepsilon \overline \theta r - \overline\varrho\overline\theta}{\varepsilon\overline\theta}\right\rangle {\bf U}\cdot\nabla F\ {\rm d}x{\rm d}t\\
+\frac 1{\overline\theta}\int_0^\tau \int_\Omega \overline \varrho F{\bf U} \cdot \nabla F + \omega(\varepsilon,\eta),
\end{multline}
where the first term on the right-hand side tends to $0$ due to \eqref{deal.pressure.term} and the second term is zero because of \eqref{Euler.boussinesq}$_1$.

Since the total energy of the acoustic equation is conserved,
\begin{multline}
\int_0^\tau \int_\Omega \frac 12 \overline\varrho \partial_t |\nabla\Phi|^2\ {\rm d}x{\rm d}t\\
 + \frac{c_v}{c_v+1}\int_0^\tau \int_\Omega\Big(\overline\theta R_{\varepsilon,\eta} + \overline\varrho T_{\varepsilon,\eta} - \overline\varrho F \Big)\partial_t \left(\frac{R_{\varepsilon,\eta}}{\overline\varrho} + \frac{T_{\varepsilon,\eta}}{\overline\theta} - \frac{F}{\overline\theta}\right)\ {\rm d}x{\rm d}t=0
\end{multline}
 and, as a result, \eqref{mezikrok.3} yields
\begin{multline}
\int_\Omega \langle U^\varepsilon_{\tau,x};\mathcal E^\varepsilon(\varrho,{\bf m}, p|\tilde \varrho_{\varepsilon,\eta},\tilde{\bf U}_{\varepsilon,\eta},\tilde \theta_{\varepsilon,\eta})\rangle\ {\rm d}x + \mathcal D^\varepsilon(\tau)\\
\leq \int_0^\tau \int_\Omega \langle U^\varepsilon_{t,x};\mathcal E^\varepsilon(\varrho,{\bf m}, p|\tilde \varrho_{\varepsilon,\eta},\tilde{\bf U}_{\varepsilon,\eta},\tilde \theta_{\varepsilon,\eta})\rangle\ + \int_0^\tau \mathcal D^\varepsilon(t)\ {\rm d}t+\omega(\varepsilon,\eta)
\end{multline}
for a.a. $\tau \in (0,T)$.
The Gronwall inequality concludes the proof.
%\newpage

\end{document}